# On some recent aspects of stochastic control and their applications*

Huyên Pham[†]

*Laboratoire de Probabilités et Modèles Aléatoires, CNRS, UMR 7599
Université Paris 7, e-mail:* `pham@math.jussieu.fr`*, and CREST*

**Abstract:** This paper is a survey on some recent aspects and developments in stochastic control. We discuss the two main historical approaches, Bellman's optimality principle and Pontryagin's maximum principle, and their modern exposition with viscosity solutions and backward stochastic differential equations. Some original proofs are presented in a unifying context including degenerate singular control problems. We emphasize key results on characterization of optimal control for diffusion processes, with a view towards applications. Some examples in finance are detailed with their explicit solutions. We also discuss numerical issues and open questions.

---

*This is an original survey paper
[†]I would like to thank Philip Protter for his encouragement to write this paper.

**Contents**

## 1. Introduction

Stochastic control is the study of dynamical systems subject to random perturbations and which can be controlled in order to optimize some performance criterion. Historically handled with Bellman's and Pontryagin's optimality principles, the research on control theory considerably developed over these last years, inspired in particular by problems emerging from mathematical finance. The global approach for studying stochastic control problems by the Bellman dynamic programming principle has now its suitable framework with viscosity solutions concept : this allows to go beyond the classical verification Bellman approach for studying degenerate singular control problems arising typically in finance. On the other hand, the stochastic maximum principle finds a modern presentation with backward stochastic differential equations (BSDE), which itself led to a very important strand of research. Viscosity solutions and BSDE have a strong relation through their PDE representation, and stochastic control theory is a place where probabilists and PDE mathematicians meet together. In this survey paper, we give an overview of some of these modern developments of stochastic control by focusing mainly on controlled diffusion processes. Our chief concern is to derive some tractable characterization of the value function and optimal control. We do not discuss the theoretical existence problem of optimal controls, which was largely studied in the literature, and we refer to [16] for a very recent overview on this subject. We include examples of typical applications with detailed explicit solutions.

The paper is organized as follows. The next section formulates the problem, discuss the classical approaches and their limitations. We also present some other stochastic control problems of interest and their possible developments. In Section 3, we present the Bellman dynamic programming approach with viscosity solutions theory. Some original proofs are exposed and we indicate for some applications in finance how examples may be explicitly solved. Section 4 describes the modern exposition of stochastic control and maximum principle by means of BSDE. An application to the stochastic linear quadratic control problem and mean-variance hedging is explicited. In both sections, we give the key results on the characterization of optimality. In Section 5, we discuss numerical issues and conclude in Section 6 with some possible developments and open questions.

## 2. The problem and discussion on methodology

### *2.1. Problem formulation*

A standard control diffusion problem on finite horizon is formulated as follows. Let $(\Omega, \mathcal{F}, P)$ be a probability space, $T > 0$ a finite time, $\mathbb{F} = (\mathcal{F}_t)_{0 \leq t \leq T}$ a filtration satisfying the usual conditions, and $W$ a $d$-dimensional Brownian motion defined on the filtered probability space $(\Omega, \mathcal{F}, \mathbb{F}, P)$. We consider the controlled state process $X$ valued in $\mathbb{R}^n$ and satisfying :

$$dX_s = b(s, X_s, \alpha_s)ds + \sigma(s, X_s, \alpha_s)dW_s. \tag{2.1}$$

The control $\alpha = (\alpha_s)_{0 \leq s \leq T}$ is a progressively measurable process valued in the control set $A$, a subset of $\mathbb{R}^m$. The Borelian functions $b$, $\sigma$ on $[0,T] \times \mathbb{R}^n \times A$ satisfy the usual conditions in order to ensure the existence of a strong solution to (2.1). This is typically satisfied when $b$ and $\sigma$ satisfy a Lipschitz condition on $(t,x)$ uniformly in $a$, and $\alpha$ satisfies a square integrability condition. We denote by $\mathcal{A}$ the set of control processes $\alpha$. Given $(t,x) \in [0,T] \times \mathbb{R}^n$ and $\alpha \in \mathcal{A}$, we then denote by $\{X_s^{t,x}, t \leq s \leq T\}$, the unique strong solution to (2.1) starting from $x$ at time $s = t$. As usual, we omit the dependance of $X$ on $\alpha$ to lighten notations.

We are then given two Borelian real-valued functions $f$ and $g$ respectively defined on $[0,T] \times \mathbb{R}^n \times A$ and $\mathbb{R}^n$ and we define the gain function :

$$J(t,x,\alpha) = E\left[\int_t^T f(s, X_s^{t,x}, \alpha_s)ds + g(X_T^{t,x})\right],$$

for all $(t,x) \in [0,T] \times \mathbb{R}^n$ and $\alpha \in \mathcal{A}$. Of course, we have to impose integrability conditions on $f$ and $g$ in order for the above expectation to be well-defined, e.g. a lower boundedness or linear growth condition. Now, since our objective is to maximize this gain function, we introduce the so-called value function :

$$v(t,x) = \sup_{\alpha \in \mathcal{A}} J(t,x,\alpha). \tag{2.2}$$

For an initial state $(t,x)$, we say that $\hat{\alpha} \in \mathcal{A}$ is an optimal control if $v(t,x) = J(t,x,\hat{\alpha})$.

**Remark 2.1** We focus mainly in this survey paper on finite horizon control problems. The infinite horizon version of control problem (2.2) is formulated as

$$\sup_{\alpha \in \mathcal{A}} E\left[\int_0^\infty e^{-\rho t} f(X_s, \alpha_s)ds\right],$$

where $\rho > 0$ is a positive discount factor, enough, to ensure finiteness of the associated value function. Such problems are studied similarly by Bellman's optimality principle as for the finite horizon case. There is also no additional difficulties for considering more general discount factors : $\rho = \rho(t,x,a)$, both in the finite or infinite horizon case. The formulation (2.2) captures some fundamental structure and properties of optimal stochastic control, but there are of course many other types of control problems that are also important both from a theoretical and applied viewpoint. We shall list some of them and present possible developments later.

### 2.2. Bellman's optimality principle

Bellman's optimality principle, initiated by Bellman [8] and also called the dynamic programming principle (DPP), is a fundamental principle in control theory : it formally means that if one has followed an optimal control decision until some arbitrary observation time, say $\theta$, then, given this information, it remains optimal to use it after $\theta$. In the context of controlled diffusion described above, the DPP is mathematically stated as follows :

(DP1) For all $\alpha \in \mathcal{A}$ and $\theta \in \mathcal{T}_{t,T}$, set of stopping times valued in $[t, T]$ :

$$v(t,x) \geq E\left[\int_t^\theta f(s, X_s^{t,x}, \alpha_s)ds + v(\theta, X_\theta^{t,x})\right]. \quad (2.3)$$

(DP2) For all $\varepsilon > 0$, there exists $\hat{\alpha}^\varepsilon \in \mathcal{A}$ s.t. for all $\theta \in \mathcal{T}_{t,T}$ :

$$v(t,x) - \varepsilon \leq E\left[\int_t^\theta f(s, X_s^{t,x}, \hat{\alpha}_s^\varepsilon)ds + v(\theta, X_\theta^{t,x})\right]. \quad (2.4)$$

Notice that this is a stronger version than the traditional version of dynamic programming principle, which is written as : for all $\theta \in \mathcal{T}_{t,T}$,

$$\text{(DPP)} \qquad v(t,x) = \sup_{\alpha \in \mathcal{A}} E\left[\int_t^\theta f(s, X_s^{t,x}, \alpha_s)ds + v(\theta, X_\theta^{t,x})\right].$$

Although the DPP has a clear intuitive meaning, its rigorous proof is technical and has been studied by several authors and by different methods. We mention among them [56], [68], [61], [15], [41] or [89]. However, it is rather difficult to find a precise reference with a self-contained and complete proof of the DPP in the above form DP1 and DP2. Formal arguments are usually given but the technical part is often omitted or does not cover exactly the considered model. Actually, it is frequent to consider originally controlled diffusion model with strong solutions of SDEs, as formulated above. Then, in order to apply the dynamic programming technique, one needs to consider a family of optimal control problems with different initial times and states along a given state trajectory. Hence, if $X$ is a state trajectory starting from $x$ at time 0 in a given filtered probability space $(\Omega, \mathcal{F}, \mathbb{F}, P)$, then $X_t$ is a random variable, for any $t > 0$, in the original probability space $(\Omega, \mathcal{F}, P)$. However, since a control is $\mathbb{F}$-adapted, the controller knows about all the relevant past information of the system up to time $t$, and in particular about $X_t$. This means that $X_t$ is not uncertain for the controller at time $t$, or mathematically, almost surely deterministic under a different probability measure $P(.|\mathcal{F}_t)$. Thus, for using dynamic programming, we are naturally led to vary the probability spaces and so to consider the weak formulation of the stochastic control problem, for which one shows the dynamic programming principle. However, it is not clearly stated in the literature that the value functions of stochastic control problems under weak and strong formulations are equal. For this paper, we shall assume the validity of the DPP in the above form and mainly focus here on the implications for the characterization of the value function through the so-called Hamilton-Jacobi-Bellman (HJB) or dynamic programming equation. We end this discussion by mentioning that Bellman's optimality principle goes beyond the framework of controlled diffusions or even Markov processes, and may also be applied for controlled semimartingales, see [31].

*2.2.1. The Hamilton-Jacobi-Bellman (in)equation*

The HJB equation is the infinitesimal version of the dynamic programming principle : it describes the local behavior of the value function $v(t,x)$ when $\theta$ is sent to $t$ in DPP. It is formally derived by assuming that the value function is $C^2$, applying Itô's formula to $v(s, X_s^{t,x})$ between $s = t$ and $s = t + h$, and then sending $h$ to zero into DPP. The classical HJB equation associated to the stochastic control problem (2.2) is :

$$-\frac{\partial v}{\partial t}(t,x) - \sup_{a \in A}[\mathcal{L}^a v(t,x) + f(t,x,a)] \;=\; 0, \quad \text{on } [0,T) \times \mathbb{R}^n, \quad (2.5)$$

where $\mathcal{L}^a$ is the second-order infinitesimal generator associated to the diffusion $X$ with constant control $a$ :

$$\mathcal{L}^a v \;=\; b(x,a).D_x v + \frac{1}{2}\text{tr}\left(\sigma(x,a)\sigma'(x,a)D_x^2 v\right).$$

This partial differential equation (PDE) is often written also as :

$$-\frac{\partial v}{\partial t}(t,x) - H(t,x,D_x v(t,x), D_x^2 v(t,x)) \;=\; 0, \quad \forall (t,x) \in [0,T) \times \mathbb{R}^n, \quad (2.6)$$

where for $(t,x,p,M) \in [0,T] \times \mathbb{R}^n \times \mathbb{R}^n \times \mathcal{S}_n$ ($\mathcal{S}_n$ is the set of symmetric $n \times n$ matrices) :

$$H(t,x,p,M) \;=\; \sup_{a \in A}\left[b(x,a).p + \frac{1}{2}\text{tr}\left(\sigma\sigma'(x,a)M\right) + f(t,x,a)\right]. \quad (2.7)$$

The function $H$ is sometimes called Hamiltonian of the associated control problem, and the PDE (2.5) or (2.6) is the dynamic programming or HJB equation. There is also an *a priori* terminal condition :

$$v(T,x) \;=\; g(x), \quad \forall x \in \mathbb{R}^n, \quad (2.8)$$

which results from the very definition of the value function $v$. We followed the usual PDE convention by writing a minus sign in front of the second order differential operator in the left side of (2.5) or (2.6), which ensures that it satisfies the ellipticity condition, see (3.3).

**Remark 2.2** The statement (2.6) of the HJB equation (in addition to the regularity problem that will be handled with viscosity solutions) requires naturally the finiteness of the Hamiltonian $H$, which is typically satisfied when the set of controls $A$ is bounded. Actually, when $A$ is unbounded, it may happen that $H(t,x,p,M)$ defined in (2.7) takes the value $\infty$ in some domain of $(t,x,p,M)$. More precisely, assume there exists a continuous function $G(t,x,p,M)$ on $[0,T] \times \mathbb{R}^n \times \mathbb{R}^n \times \mathcal{S}_n$ s.t.

$$H(t,x,p,M) \;<\; \infty \quad \Longleftrightarrow \quad G(t,x,p,M) \;\geq\; 0.$$

Then, from (2.6), we must have :

$$G(t, x, D_x v(t,x), D_x^2 v(t,x)) \geq 0, \quad (2.9)$$

$$\text{and} \quad -\frac{\partial v}{\partial t}(t,x) - H(t, x, D_x v(t,x), D_x^2 v(t,x)) \geq 0 \quad (2.10)$$

Moreover, if inequality (2.9) is strict, then inequality (2.10) is an equality. Hence, we formally get the HJB variational inequality associated to the dynamic programming principle :

$$\min\left[-\frac{\partial v}{\partial t}(t,x) - H(t, x, D_x v(t,x), D_x^2 v(t,x)),\right.$$
$$\left. G(t, x, D_x v(t,x), D_x^2 v(t,x))\right] = 0.$$

In this case, we say that the control problem is singular, in contrast with the regular case of HJB equation (2.6). A typical example of a singular problem occurs when the control influences linearly the dynamics of the system and the gain function. To fix the idea, consider the one-dimensional case $n = 1$, $A = \mathbb{R}_+$, and

$$b(x,a) = \hat{b}(x) + a, \quad \sigma(x,a) = \hat{\sigma}(x), \quad f(t,x,a) = \hat{f}(t,x) - \lambda a,$$

for some $\lambda \in \mathbb{R}$. Then,

$$H(t,x,p,M) = \begin{cases} \hat{b}(x)p + \frac{1}{2}\hat{\sigma}(x)^2 M + \hat{f}(t,x) & \text{if } -p + \lambda \geq 0 \\ \infty & \text{if } -p + \lambda < 0. \end{cases}$$

The HJB variational inequality is then written as :

$$\min\left[-\frac{\partial v}{\partial t}(t,x) - \hat{b}(x)\frac{\partial v}{\partial x}(t,x) - \frac{1}{2}\hat{\sigma}(x)^2\frac{\partial^2 v}{\partial x^2}(t,x) , -\frac{\partial v}{\partial x}(t,x) + \lambda\right] = 0.$$

We shall give in the next section another example of singular control arising in finance and where the control is in the diffusion term.

In a singular control problem, the value function is in general discontinuous in $T$ so that (2.8) is not the relevant terminal condition. We shall show how to derive rigorously the HJB equation (or variational inequality) with the concept of viscosity solution to handle the lack of *a priori* regularity of the value function, and also to determine the correct terminal condition.

*2.2.2. The classical verification approach*

The classical verification approach consists in finding a smooth solution to the HJB equation, and to check that this candidate, under suitable sufficient conditions, coincides with the value function. This result is usually called a verification theorem and provides as a byproduct an optimal control. It relies mainly on Itô's formula. The assertions of a verification theorem may slightly vary from problem to problem, depending on the required sufficient technical conditions. These conditions should actually be adapted to the context of the considered problem. In the above context, a verification theorem is roughly stated as follows :

VERIFICATION THEOREM. Let $w$ be a $C^{1,2}$ function on $[0,T) \times \mathbb{R}^n$ and continuous in $T$, with suitable growth condition. Suppose that for all $(t,x) \in [0,T) \times \mathbb{R}^n$, there exists $\hat{\alpha}(t,x)$ mesurable, valued in $A$ s.t. $w$ solves the HJB equation :

$$\begin{aligned} 0 &= -\frac{\partial w}{\partial t}(t,x) - \sup_{a \in A}\left[\mathcal{L}^a w(t,x) + f(t,x,a)\right] \\ &= -\frac{\partial w}{\partial t}(t,x) - \mathcal{L}^{\hat{\alpha}(t,x)} w(t,x) - f(t,x,\hat{\alpha}(t,x)) \quad \text{on} \quad [0,T) \times \mathbb{R}^n, \end{aligned}$$

together with the terminal condition

$$w(T,.) = g \quad \text{on} \quad \mathbb{R}^n,$$

and the S.D.E. :

$$dX_s = b(X_s, \hat{\alpha}(s,X_s))ds + \sigma(X_s, \hat{\alpha}(s,X_s))dW_s$$

admits a unique solution, denoted $\hat{X}_s^{t,x}$, given an initial condition $X_t = x$. Then, $w = v$ and $\{\hat{\alpha}(s, \hat{X}_s^{t,x}) \ t \leq s \leq T\}$ is an optimal control for $v(t,x)$.

A proof of this verification theorem may be found in any textbook on stochastic control, see e.g. [39], [56], [41], [89] or [81]. The first and most famous application in finance of this verification theorem for stochastic control problem is Merton's portfolio selection problem. This is the situation where an investor may decide at any time over a finite horizon $T$ to invest a proportion $\alpha$ valued in $A = \mathbb{R}$ of his wealth $X$ in a risky stock of constants rate of return $\mu$ and volatility $\sigma$ and the rest of proportion $1 - \alpha$ in a bank account of constant interest $r$. His wealth controlled process is then governed by :

$$dX_s = X_s \left(r + (\mu - r)\alpha_s\right)ds + X_s \sigma \alpha_s dW_s,$$

and the objective of the investor is given by the value function :

$$v(t,x) = \sup_{\alpha \in \mathcal{A}} E\left[U(X_T^{t,x})\right], \quad (t,x) \in [0,T] \times \mathbb{R}_+,$$

where $U$ is a utility function, i.e. a concave and increasing function on $\mathbb{R}_+$. For the popular specific choice of the power utility function $U(x) = x^p$, with $p < 1$, it is possible to find an explicit (smooth) solution to the associated HJB equation with the terminal condition $v(T,.) = U$, namely :

$$v(t,x) = \exp\left(\rho(T-t)\right)x^p,$$

with $\rho = \frac{(\mu-r)^2}{2\sigma^2}\frac{p}{1-p} + rp$. Moreover, the optimal control is constant and given by :

$$\hat{\alpha} = \frac{\mu - r}{\sigma^2(1-p)}$$

The key point in the explicit resolution of the Merton problem is that the value function $v$ may be separated into a function of $t$ and of $x$: $v(t,x) = \varphi(t)x^p$. With this transformation and substitution into the HJB equation, it turns out that $\varphi$ is the solution of an ordinary differential equation with terminal condition $\varphi(T) = 1$, which is explicitly solved. Other applications of verification theorems to stochastic control problems in finance are developed in the recent textbook [70]. There are also important examples of applications in engineering with the stochastic linear regulator, which will be studied later by the maximum principle.

Classical verification theorems allow to solve examples of control problems where one can find, or at least there exists, a smooth solution to the associated HJB equation. They apply successfully for control problems where the diffusion term does not depend on the control and is uniformly elliptic, since in this case the HJB equation is semilinear in the sense that $D_x^2 v$ appears linearly, and so classical existence results for smooth solutions exist, see [60]. They also apply for some specific models with control on the diffusion term, typically Merton's portfolio selection problem as described above, and more generally to extensions of Merton's model with stochastic volatility, since in this case, the HJB equation may be reduced after a suitable transformation to a semilinear equation. This last point is developed in [91] and [77]. However, in the general case of nondegeneracy of the diffusion term and in a number of applications, the value function might be not smooth or it is not possible to obtain *a priori* the required regularity. Moreover, for singular control problems, the value function is in general not continuous at the terminal date $T$, so that the right terminal condition is not given by $g(x)$, i.e. $\lim_{t \nearrow T} v(t,x) \neq g(x)$. Then, the classical verification approach does not work and we need to relax the notion of solution to the HJB equation. It turns out that the suitable class of weak solutions is the one of viscosity solutions, not only for dealing with the rigorous derivation of the HJB equation, but also for determining the relevant terminal condition. Moreover, viscosity solutions theory provides a general verification approach, which allows to go beyond the case of classical verification. We shall give some applications in finance of the viscosity solutions approach where we can explicit non smooth solutions to singular control problems.

The notion of viscosity solutions was introduced by P.L. Lions [61] for second order equations. They provide by now a well-established method for dealing with stochastic control problems, in particular in finance, see e.g. [90], [41], [89] or [81]. We review some of these aspects and their applications in the next section.

### 2.3. Pontryagin's maximum principle

A classical approach for optimization and control problems is to derive necessary conditions satisfied by an optimal solution. For example, the Kuhn-Tucker condition is a well-known necessary condition of optimality in the finite-dimensional case. For an optimal control problem as described above, and which may be viewed as an infinite-dimensional optimization problem, the argument is to use an appropriate calculus of variations on the gain function $J(t, x, \alpha)$ with respect to the control variable in order to derive a necessary condition of optimality. The maximum principle, initiated by Pontryagin in the 1960s, states that an optimal state trajectory must solve a Hamilton system together with a maximum condition of a function called a generalized Hamilton. In principle, solve a Hamilton should be easier than solving the original control problem.

The original version of Pontryagin's maximum principle was derived for deterministic problems. As in classical calculus of variation, the basic idea of is to perturb an optimal control and to use some sort of Taylor expansion of the state trajectory and objective functional around the optimal control. By sending the perturbation to zero, one obtains some inequality, and by duality, the maximum principle is expressed in terms of an adjoint variable (Lagrange multiplier in the finite-dimensional case). The stochastic control case was extensively studied in the 1970s by Bismut, Kushner, Bensoussan or Haussmann. However, at that time, the results were essentially obtained under the condition that there is no control on the diffusion coefficient. For example, Haussmann investigated maximum principle by Girsanov's transformation and this limitation explains why this approach does not work with control-dependent and degenerate diffusion coefficients. The main difficulty when facing a general controlled diffusion is that the Itô integral term is not of the same order as the Lebesgue term and thus the first-order variation method fails. This difficulty was overcomed by Peng [75], who studied the second-order term in the Taylor expansion of the perturbation method arising from the Itô integral. He then obtained a maximum principle for possibly degenerate and control-dependent diffusion, which involves in addition to the first-order adjoint variable, a second-order adjoint variable.

In order to make applicable the maximum principle, one needs some explicit description of the adjoint variables. These variables obtained originally by duality in functional analysis may be represented by Riesz representation of a certain functional. By completing with martingale representation in stochastic analysis, the adjoint variables are then described by what is called today backward stochastic differential equations (BSDE). Actually, it was the study of the maximum principle in the stochastic control case that motivated Peng for a general formulation of BSDE, which in turn generated an important area of research. We shall state in Section 4 the mathematical formulation of the sufficiency of the maximum principle and relate it to BSDE and its extensions.

### *2.4. Other control problems*

We present in this paragraph some other control problems, which we do not study in detail here, but that are also of significant theoretical and practical interest. We also emphasize some present developments.

#### *2.4.1. Random horizon*

In problem formulation (2.2), the time horizon is fixed until a deterministic terminal time $T$. In some real applications, the time horizon may be random, and in the context of controlled diffusion described in paragraph 2.1, the control problem is formulated as :

$$\sup_{\alpha \in \mathcal{A}} E \left[ \int_0^\tau f(s, X_s, \alpha_s) ds + g(X_\tau) \right], \tag{2.11}$$

where $\tau$ is a finite random time. In standard cases, the terminal time $\tau$ is a stopping time at which the state process exits from a certain relevant domain. For example, in a reinsurance model, the state process $X$ is the reserve of a company that may control it by reinsuring a proportion $1 - \alpha$ of premiums to another company. The terminal time $\tau$ is then the bankruptcy time of the company defined as $\tau = \inf\{t \geq 0 \ : \ X_t \leq 0\}$. More generally, given some open set $\mathcal{O}$ of $\mathbb{R}^n$,

$$\tau \ = \ \inf\{t \geq 0 \ : \ X_t \notin \mathcal{O}\} \wedge T$$

(which depends on the control). In this case, the control problem (2.11) leads via the dynamic programming approach to a Dirichlet boundary-value problem. Another case of interest concerns a terminal time $\tau$, which is a random time but not a stopping time in the filtration $\mathbb{F}$ with respect to which the controls are adapted. This situation occurs for example in credit risk models where $\tau$ is the default time of a firm. Under the so-called **(G)** hypothesis on filtration theory : $P[\tau \leq t | \mathcal{F}_t]$ is a nondecreasing right-continuous process, problem (2.11) may be reduced to a stochastic control problem under a fixed deterministic horizon, see [13] for a recent application in portfolio optimization model. In the general random time case, the associated control problem has been relatively lightly studied in the literature, see [17] or [92] for a utility maximization problem in finance.

### 2.4.2. *Optimal stopping*

In the models presented above, the horizon of the problem is either fixed or indirectly influenced by the control. When one has the possibility to control directly the terminal time, which is then modelled by a controlled stopping time, the associated problem is an optimal stopping time problem. In the general formulation of such models, the control is mixed, composed by a pair control/stopping time $(\alpha, \tau)$ and the functional to optimize is :

$$E \left[ \int_0^\tau f(t, X_t, \alpha_t) dt + g(X_\tau) \right].$$

The theory of optimal stopping, thoroughly studied in the seventies, has received a renewed interest with a variety of applications in economics and finance. These applications range from asset pricing (American options) to firm investment and real options. Extensions of classical optimal stopping problems deal with multiple optimal stopping with eventual changes of regimes in the state process. They were studied e.g. in [9], [87], and applied in finance in [19], [30], [45], [21] or [80].

*2.4.3. Impulse control*

In formulation of the problem in paragraph 2.1, the displacement of the state changes continuously in time in response to the control effort. However, in many real applications, this displacement may be discontinuous. For example, in insurance company models, the company distributes the dividends once or twice a year rather than continuously. In transaction costs models, the agent should not invest continuously in the stock due to the costs but only at discrete times. A similar situation occurs in a liquidity risk model, see e.g. [22]. Impulse control provides a suitable framework for modelling such situations. This may be described as follows : the controlled state diffusion process is governed by

$$dX_s \;=\; b(s,X_s)dt + \sigma(s,X_s)dW_s + d\zeta_s,$$

where the control $\zeta$ is a pure jump process. In other words, the control is given by a pair $(\tau_n, \kappa_n)_n$ where $(\tau_n)_n$ is a nondecreasing sequence of stopping times, representing the intervention times of the controller, and $(\kappa_n)_n$ is a sequence of $\mathcal{F}_{\tau_n}$-measurable random variables, representing the jump size decided by the controller at time $\tau_n$. The functional objective to optimize is in the form :

$$E\left[\int_0^T f(t,X_t,\alpha_t)dt + \sum_{\tau_n \leq T} h(X_{\tau_n},\kappa_n) + g(X_T)\right].$$

Impulse control problem is known to be associated via the dynamic programming approach to an HJB quasi-variational inequality, see [9]. For some recent applications in finance, we refer to [49] for insurance models, [55] and [69] for transaction costs models, and [63] for liquidity risk model.

*2.4.4. Ergodic control*

Some stochastic systems may exhibit over a long period a stationary behavior characterized by an invariant measure. This measure, if it does exists, is obtained by the average of the states over a long time. An ergodic control problem consists in optimizing over the long term some criterion taking into account this invariant measure.

A standard formulation resulting from the criterion presented in paragraph 2.1 is to optimize over control $\alpha$ functional of the form :

$$\limsup_{T\to+\infty} \frac{1}{T} E\left[\int_0^T f(X_t,\alpha_t)dt\right],$$

or

$$\limsup_{T\to+\infty} \frac{1}{T} \ln E\left[\exp\left(\int_0^T f(X_t,\alpha_t)dt\right)\right].$$

This last formulation is called risk-sensitive control on an infinite horizon. Ergodic and risk-sensitive control problems were studied in [51], [11] or [38]. Risk-sensitive control problems were recently applied in a financial context in [12] and [40].

Another criterion is based on the large deviations behavior of the ergodic system : $P[X_T/T] \simeq e^{-I(c)T}$, when $T$ goes to infinity, consists in maximizing over control $\alpha$ functional of the form :

$$\limsup_{T\to+\infty} \frac{1}{T} \ln P\left[\frac{X_T}{T} \geq c\right].$$

This large deviations control problem is interpreted in finance as the asymptotic version of the quantile criterion of maximizing the probability that the terminal wealth $X_T$ beats a given benchmark. This nonstandard control problem was introduced and developed recently by Pham [78], [79], see also [47]. It does not have a direct dynamic programming principle but may be reduced via a duality principle to a risk-sensitive control problem.

*2.4.5. Robust control*

In the problems formulated above, the dynamics of the control system is assumed to be known and fixed. Robust control theory is a method to measure the performance changes of a control system with changing system parameters. This is of course important in engineering systems, and it has recently been used in finance in relation with the theory of risk measures initiated by Artzner et al. [3]. Indeed, it is proved that a coherent risk measure for an uncertain payoff $X_T$ at time $T$ is represented by :

$$\rho(-X_T) = \sup_{Q\in\mathcal{Q}} E^Q[X_T],$$

where $\mathcal{Q}$ is a set of absolutely continuous probability measures with respect to the original probability $P$. More generally, one may define a risk measure in the form :

$$\rho(-X_T) = -\inf_{Q\in\mathbb{Q}} E^Q[U(-X_T)],$$

where $U$ is a concave and nondecreasing function. So, when $X$ is controlled by $\alpha$ and the problem is to optimize the risk measure $\rho(-X_T)$, one is facing a robust control problem. In this financial context, robust optimization problems were recently studied in [82] and [44].

*2.4.6. Partial observation control problem*

It is assumed so far that the controller completely observes the state system. In many real applications, he is only able to observe partially the state via other variables and there is noise in the observation system. For example in financial models, one may observe the asset price but not completely its rate of return and/or its volatility, and the portfolio investment is based only on the asset price information. We are facing a partial observation control problem. This may be formulated in a general form as follows : we have a controlled signal (unobserved) process governed by

$$dX_s = b(s, X_s, Y_s, \alpha_s)ds + \sigma(s, X_s, Y_s, \alpha_s)dW_s,$$

and an observation process

$$dY_s = \eta(s, X_s, Y_s, \alpha_s)ds + \gamma(s, X_s, Y_s, \alpha_s)dB_s,$$

where $B$ is another Brownian motion, eventually correlated with $W$. The control $\alpha$ is adapted with respect to the filtration generated by the observation $\mathbb{F}^Y = (\mathcal{F}^Y_t)$ and the functional to optimize is :

$$J(\alpha) = E\left[\int_0^T f(X_t, Y_t, \alpha_t)dt + g(X_T, Y_T)\right].$$

By introducing the filter measure-valued process

$$\Pi_t(dx) = P[X_t \in dx | \mathcal{F}^Y_t],$$

one may rewrite the functional $J(\alpha)$ in the form :

$$J(\alpha) = E\left[\int_0^T \hat{f}(\Pi_t, Y_t, \alpha_t)dt + \hat{g}(\Pi_T, Y_T)\right],$$

where we use the notation : $\hat{f}(\pi, y) = \int f(x, y)\pi(dx)$ for any finite measure $\pi$ on the signal state space, and similarly for $\hat{g}$. Since by definition, the process $(\Pi_t)$ is $(\mathcal{F}^Y_t)$-adapted, the original partial observation control problem is reformulated as a complete observation control model, with the new observable state variable defined by the filter process. The additional main difficulty is that the filter process is valued in the infinite-dimensional space of probability measures : it satisfies the Zakai stochastic partial differential equation. The dynamic programming principle or maximum principle are still applicable and the associated Bellman equation or Hamiltonian system are now in infinite dimension. For a theoretical study of optimal control under partial observation under this infinite dimensional viewpoint, we mention among others the works [36], [26], [5], [10], [62] or [93]. There are relatively few explicit calculations in the applications to finance of partial observation control models and this area should be developed in the future.

*2.4.7. Stochastic target problems*

Motivated by the superreplication problem in finance, and in particular under gamma constraints [85], Soner and Touzi introduced a new class of stochastic control problems. The state process is described by a pair $(X, Y)$ valued in $\mathbb{R}^n \times \mathbb{R}$, and controlled by a control process $\alpha$ according to :

$$dX_s = b(s, X_s, \alpha_s)ds + \sigma(s, X_s, \alpha_s)dW_s \qquad (2.12)$$
$$dY_s = \eta(s, X_s, Y_s, \alpha_s)ds + \gamma(s, X_s, Y_s, \alpha_s)dW_s. \qquad (2.13)$$

Notice that the coefficients of $X$ do not depend on $Y$. Given $\alpha \in \mathcal{A}$ and $(t, x, y) \in [0, T] \times \mathbb{R}^n \times \mathbb{R}$, $(X^{t,x}, Y^{t,x,y})$ is the unique solution to (2.12)-(2.13) with initial condition $(X_t^{t,x}, Y_t^{t,x,y}) = (x, y)$. The coefficients $b, \sigma, \eta, \gamma$ are bounded functions and satisfy usual conditions ensuring that $(X^{t,x}, Y^{t,x,y})$ is well-defined. The stochastic target problem is defined as follows. Given a real-valued measurable function $g$ on $\mathbb{R}^n$, the value function of the control problem is defined by :

$$v(t, x) = \inf \left\{ y \in \mathbb{R} : \exists \alpha \in \mathcal{A}, \ Y_T^{t,x,y} \geq g(X_T^{t,x}) \ a.s. \right\}.$$

In finance, $X$ is the price process, $Y$ is the wealth process controlled by the portfolio strategy $\alpha$, and $v(t, x)$ is the minimum capital in order to superreplicate the payoff option $g(X_T)$.

The dynamic programming principle associated to this stochastic target problem is stated as follows : for all $(t, x) \in [0, T]$, and $\theta$ stopping times in $[t, T]$, we have

$$v(t, x) = \inf \left\{ y \in \mathbb{R} : \exists \alpha \in \mathcal{A}, \ Y_\theta^{t,x,y} \geq v(\theta, X_\theta^{t,x}) \ a.s. \right\}.$$

The derivation of the associated dynamic programming equation is obtained under the following conditions. The matrix $\sigma(t, x, a)$ is assumed to be invertible and the function

$$a \mapsto \gamma(t, x, y, a)\sigma^{-1}(t, x, a)$$

is one-to-one for all $(t, x, y)$, with inverse denoted $\vartheta$, i.e.

$$\gamma(t, x, y, a)\sigma^{-1}(t, x, a) = p' \iff a = \vartheta(t, x, y, p)$$

for $p \in \mathbb{R}^n$. Here $p'$ is the transpose of $p$. Moreover, the control set $A$ is assumed to be convex and compact in $\mathbb{R}^m$, with nonempty interior. The support function of the closed convex set $A$ is denoted

$$\delta_A(\zeta) = \sup_{a \in A} a'\zeta.$$

Under these conditions, Soner and Touzi [86] proved that $v$ is a viscosity solution to the dynamic programming equation :

$$\min \left\{ -\frac{\partial v}{\partial t}(t, x) - \mathcal{L}^{a_0(t,x)} v(t, x) + \eta(t, x, v(t, x), a_0(t, x)) \right., $$
$$\left. G(t, x, v(t, x), D_x v(t, x)) \right\} = 0,$$

where $\mathcal{L}^a$ is the second order differential operator associated to the diffusion $X$, and

$$a_0(t, x) = \vartheta(t, x, v(t, x), D_x v(t, x)),$$
$$G(t, x, v(t, x), D_x v(t, x)) = \inf_{|\zeta|=1} \left[ \delta_A(\zeta) - a_0(t, x)'\zeta \right].$$

## 3. Dynamic programming and viscosity solutions

### 3.1. Definition of viscosity solutions

The notion of viscosity solutions provides a powerful tool for handling stochastic control problems. The theory of viscosity solutions goes beyond the context of HJB equations and is a concept of general weak solutions for partial differential equations. We refer to the user's guide of Crandall, Ishii and Lions [23] for an overview of this theory. Here, we simply recall the definition and some basic properties required for our purpose.

Let us consider parabolic nonlinear partial differential equations of second order :

$$F(t,x,w(t,x),\frac{\partial w}{\partial t}(t,x),D_x w(x),D^2_{xx}w(x)) \; = \; 0, \quad (t,x) \in [0,T) \times \mathcal{O} \quad (3.1)$$

where $\mathcal{O}$ is an open set of $\mathbb{R}^n$, and $F$ is a continuous function on $[0,T] \times \mathcal{O} \times \mathbb{R} \times \mathbb{R} \times \mathbb{R}^n \times \mathcal{S}_n$. The function $F$ is assumed to satisfy the ellipticity and parabolicity conditions :

$$M \leq \widehat{M} \quad \Longrightarrow \quad F(t,x,r,p_t,p,M) \geq F(t,x,r,p_t,p,\widehat{M}) \qquad (3.2)$$

$$p_t \leq \hat{p}_t \quad \Longrightarrow \quad F(t,x,r,p_t,p,M) \geq F(t,x,r,\hat{p}_t,p,M), \qquad (3.3)$$

for all $t \in [0,T)$, $x \in \mathcal{O}_n$, $r \in \mathbb{R}$, $p_t, \hat{p}_t \in \mathbb{R}$, $p \in \mathbb{R}$ and $M, \hat{M} \in \mathcal{S}_n$. The last condition (3.3) means that we are dealing with forward PDE, i.e. (3.1) holds for time $t < T$, and the terminal condition is for $t = T$. This is in accordance with the control problem and HJB equation formulated in Section 2.

Since it is not always easy to have *a priori* the continuity of the value function (actually, there may be even discontinuity at terminal time $T$), we work with the notion of discontinuous viscosity solutions. We then introduce, for a locally bounded function $w$ on $[0,T] \times \mathcal{O}$, its lower semicontinuous envelope $w_*$ and upper semicontinuous envelope $w^*$, i.e. :

$$w_*(t,x) \; = \; \liminf_{(t',x')\to(t,x)} w(t',x') \quad \text{and} \quad w^*(x) \; = \; \limsup_{(t',x')\to(t,x)} w(t',x').$$

**Definition 3.1** *Let $w$ be a locally bounded function on $\mathcal{O}$.*
*(i) $w$ is a viscosity subsolution (resp. supersolution) of (3.1) if :*

$$F(\bar{t},\bar{x},\varphi(\bar{t},\bar{x}),\frac{\partial \varphi}{\partial t}(\bar{t},\bar{x}),D_x\varphi(\bar{t},\bar{x}),D^2_x\varphi(\bar{t},\bar{x})) \; \leq \; \text{(resp. } \geq \text{)} \; 0, \quad (3.4)$$

*for any $(\bar{t},\bar{x}) \in [0,T) \times \mathcal{O}$ and smooth test function $\varphi \in C^2([0,T] \times \mathcal{O})$ s.t. $(\bar{t},\bar{x})$ is a maximum (resp. minimum) of $w^* - \varphi$ (resp. $w_* - \varphi$) with $0 = (w^* - \varphi)(\bar{t},\bar{x})$ (resp. $(w_* - \varphi)(\bar{t},\bar{x})$).*
*(ii) $w$ is a viscosity solution of (3.1) if it is a viscosity subsolution and supersolution.*

**Remark 3.1 1.** Viscosity solutions extend the notion of classical solutions : a $C^{1,2}$ function on $[0,T] \times \mathcal{O}$ is a supersolution (resp. a subsolution) of (3.1) in the classical sense iff it is a viscosity supersolution (resp. viscosity subsolution).
**2.** The above definition is unchanged if the maximum/minimum is strict and/or local. Notice that since we are considering forward PDEs, a local extremum point $(\bar{t}, \bar{x})$ means with respect to a neighborhood of the form $[\bar{t}, \bar{t}+h) \times B_\eta(\bar{x})$, with $h, \eta > 0$. Here $B_\eta(x)$ denotes the open ball of radius $\eta$ and center $x$ and $\bar{B}_\eta(x)$ its closure.
**3.** Without loss of generality, by translating the test function by a constant, we may relax the condition that $0 = (w^* - \varphi)(\bar{t}, \bar{x})$ (resp. $(w_* - \varphi)(\bar{t}, \bar{x})$). We have then to replace in (3.4) $\varphi(\bar{t}, \bar{x})$ by $w^*(\bar{t}, \bar{x})$ (resp. $w_*(\bar{t}, \bar{x})$).
**4.** We define similarly viscosity solutions for elliptic PDEs, i.e. without time variable $t$.

### 3.2. Viscosity characterization

We come back to the framework of controlled diffusions of Section 1, and we state that the value function is a viscosity solution to the associated HJB equation. We also determine the relevant terminal condition. We present a unifying result for taking into account the possible singularity of the Hamiltonian $H$ when the control set $A$ is unbounded, see Remark 2.2. We then introduce

$$\text{dom}(H) = \{(t,x,p,M) \in [0,T] \times \mathbb{R}^n \times \mathbb{R}^n \times \mathcal{S}_n : H(t,x,p,M) < \infty\},$$

and we shall make the assumption :

$$\begin{aligned} &H \text{ is continuous on int}(\text{dom}(H)) \\ &\text{and there exists } G : [0,T] \times \mathbb{R}^n \times \mathbb{R}^n \times \mathcal{S}_n \\ &\text{nondecreasing in its last argument and continuous s.t. :} \\ &(t,x,p,M) \in \text{dom}(H) \iff G(t,x,p,M) \geq 0 \end{aligned} \quad (3.5)$$

The nondecreasing condition of $G$ in its last argument is only required here to ensure that the PDE in (3.6) satisfies the ellipticity condition (3.2).

In the sequel, we assume that the value function $v$ is locally bounded on $[0,T] \times \mathbb{R}^n$ : this is a condition usually easy to check, and satisfied typically when $v$ inherits from $f$ and $g$ a linear growth condition. Our general viscosity property for the value function is stated in the following theorem :

**Theorem 3.1** *Assume that $f(.,.,a)$ is continuous for all $a \in A$ and (3.5) holds. Then, $v$ is a viscosity solution of the HJB variational inequality :*

$$\min\left\{-\frac{\partial v}{\partial t}(t,x) - H(t,x,D_x v(t,x), D_x^2 v(t,x)),\right.$$
$$\left. G(t,x,D_x v(t,x), D_x^2 v(t,x))\right\} = 0, \quad (t,x) \in [0,T] \times \mathbb{R}^n. \,(3.6)$$

**Remark 3.2** In the regular case, i.e. when the Hamiltonian $H$ is finite on the whole state domain (this occurs typically when the control set $A$ is compact), the condition (3.5) is satisfied for any choice of positive continuous function, e.g. a positive constant. In this case, the HJB variational inequality (3.6) is reduced to the regular HJB equation :

$$-\frac{\partial v}{\partial t}(t,x) - H(t,x,D_xv(t,x),D_x^2v(t,x)) \;=\; 0, \quad (t,x)\in[0,T)\times\mathbb{R}^n$$

that the value function satisfies in the viscosity sense. Hence, Theorem 3.1 states a general viscosity property including both the regular and singular case.

**Remark 3.3** We do not address here the important uniqueness problem associated to the HJB (in)equation. We refer to [23] for some general uniqueness results. In most cases, there is a (strong) comparison principle for this HJB PDE, which states that any uppersemicontinuous subsolution is not greater than a lowersemicontinuous supersolution. This implies $v^* \leq v_*$ and therefore $v^* = v_*$ since the other inequality is always true by definition. The consequence of this is the continuity of the value function $v$ on $[0,T)\times\mathbb{R}^n$. Hence, we notice that the viscosity property allows to derive as a byproduct the continuity of the value function, which may not always be easily proved by a direct probabilistic argument.

It is well known that to a parabolic PDE is associated a terminal condition, in order to get a uniqueness result. We then need to determine the terminal data for the value function. By the very definition of the value function, we have

$$v(T,x) \;=\; g(x), \quad x\in\mathbb{R}^d. \tag{3.7}$$

However, in several aplications of stochastic control problems, the value function may be discontinuous at $T$, see e.g. [20] and [24]. In this case, (3.7) is not the relevant terminal condition associated to the HJB equation in order to characterize the value function : we need actually to determine $v(T^-,x) := \lim_{t\nearrow T} v(t,x)$ if it exists.

To this end, we introduce

$$\underline{v}(x) \;=\; \liminf_{t\nearrow T, x'\to x} v(t,x') \quad \text{and} \quad \bar{v}(x) \;=\; \limsup_{t\nearrow T, x'\to x} v(t,x')$$

Notice that by definition, $\underline{v} \leq \bar{v}$, $\underline{v}$ is lowersemicontinuous, and $\bar{v}$ is uppersemicontinuous.

Our characterization result for the terminal condition is stated as follows.

**Theorem 3.2** *Assume that $f$ and $g$ are lower-bounded or satisfy a linear growth condition, and (3.5) holds.*
*(i) Suppose that $g$ is lowersemicontinuous. Then $\underline{v}$ is a viscosity supersolution of*

$$\min\left[\underline{v}(x) - g(x)\,,\; G(T,x,D\underline{v}(x),D^2\underline{v}(x))\right] \;=\; 0, \quad on\;\; \mathbb{R}^n. \tag{3.8}$$

*(ii) Suppose that g is uppersemicontinuous. Then $\bar{v}$ is a viscosity subsolution of*

$$\min \left[ \bar{v}(x) - g(x) \, , \, G(T, x, D\bar{v}(x), D^2\bar{v}(x)) \right] \;\; = \;\; 0, \quad on \;\; \mathbb{R}^n. \qquad (3.9)$$

**Remark 3.4** In most cases, there is a comparison principle for the PDE arising in the above theorem, meaning that a subsolution is not greater than a supersolution. Therefore, under the conditions of Theorem 3.2, we have $\bar{v} \leq \underline{v}$ and so $\bar{v} = \underline{v}$. This means $\hat{v} := v(T^-, .)$ exists, equal to $\underline{v} = \bar{v}$ and is a viscosity solution to :

$$\min \left[ \hat{v}(x) - g(x) \, , \, G(T, x, D\hat{v}(x), D^2\hat{v}(x)) \right] \;\; = \;\; 0, \quad \text{on } \mathbb{R}^n. \qquad (3.10)$$

Denote, by $\hat{g}$ the upper $G$-envelope of $g$, defined as the smallest function above $g$ and viscosity supersolution to :

$$G(T, x, D\hat{g}(x), D^2\hat{g}(x)) \;\; = \;\; 0, \quad \text{on } \mathbb{R}^n, \qquad (3.11)$$

when it exists and is finite. Such a function may be calculated in a number of examples, see e.g. paragraph 3.3.2. Since $\hat{v}$ is a viscosity supersolution to (3.10), it is greater than $g$ and is a viscosity supersolution to the same PDE as $\hat{g}$. Hence, by definition of $\hat{g}$, we have $\hat{v} \geq \hat{g}$. On the other hand, $\hat{g}$ is a viscosity supersolution to the PDE (3.10), and so by a comparison principle, the subsolution $\hat{v}$ of (3.10) is not greater than $\hat{g}$. We have then determined explicitly the terminal data :

$$v(T^-, x) \;\; = \;\; \hat{v}(x) \;\; = \;\; \hat{g}(x).$$

Recall that in the regular case, we may take for $G$ a positive constant function, so that obviously $\hat{g} = g$. Therefore, in this case, $v$ is continuous in $T$ and $v(T^-, x) = v(T, x) = g(x)$. In the singular case, $\hat{g}$ is in general different from $g$ and so $v$ is discontinuous in $T$. The effect of the singularity is to lift up, via the $G$ operator, the terminal function $g$ to $\hat{g}$.

We separate the proof of viscosity supersolution and subsolution, which are quite different. The supersolution part follows from DP1 and standard arguments in viscosity solution theory as in Lions [61]. In this part, the notion of viscosity solutions is only used to handle the lack of *a priori* regularity of $v$. The subsolution part is more delicate and should take into account the possible singular part of the Hamiltonian. The derivation is obtained from DP2 and a contraposition argument. The assertion of the subsolution part as well as the determination of the terminal condition seem to be new in this unifying context and are inspired by arguments in [86]. These arguments really use the concept of viscosity solutions even if the value function were known to be smooth. The reader who is not interested in the technical proofs of Theorems 3.1 and 3.2 can go directly to paragraph 3.3 for some applications in finance of viscosity solutions theory.

*3.2.1. Viscosity supersolution property*

**Proposition 3.1** *Assume that $f(.,.,a)$ is continuous for all $a \in A$. Then, $v$ is a viscosity supersolution of the HJB equation :*

$$-\frac{\partial v}{\partial t}(t,x) - H(t,x,D_x v(t,x), D_x^2 v(t,x)) = 0, \quad (t,x) \in [0,T) \times \mathbb{R}^n. \quad (3.12)$$

**Proof.** Let $(\bar{t}, \bar{x}) \in [0,T) \times \mathbb{R}^n$ and $\varphi \in C^2([0,T) \times \mathbb{R}^n)$ a smooth test function s.t. :

$$0 = (v_* - \varphi)(\bar{t}, \bar{x}) = \min_{(t,x) \in [0,T) \times \mathbb{R}^n} (v_* - \varphi)(t,x). \quad (3.13)$$

By definition of $v_*(\bar{t}, \bar{x})$, there exists a sequence $(t_m, x_m)$ in $[0,T) \times \mathbb{R}^n$ s.t.

$$(t_m, x_m) \to (\bar{t}, \bar{x}) \quad \text{and} \quad v(t_m, x_m) \to v_*(\bar{t}, \bar{x}),$$

when $m$ goes to infinity. By continuity of $\varphi$ and (3.13), we also have

$$\gamma_m := v(t_m, x_m) - \varphi(t_m, x_m) \to 0,$$

when $m$ goes to infinity.

Let $a$ be an arbitrary element in $A$, and $\alpha$ the constant control equal to $a$. We denote by $X_s^{t_m, x_m}$ the associated controlled process starting from $x_m$ at $t_m$. Consider $\tau_m$ the first exit time of $X^{t_m, x_m}$ from the open ball $B_\eta(x_m)$ : $\tau_m = \inf\{s \geq t_m : |X_s^{t_m, x_m} - x_m| \geq \eta\}$, with $\eta > 0$ and let $(h_m)$ a positive sequence s.t. :

$$h_m \to 0 \quad \text{and} \quad \frac{\gamma_m}{h_m} \to 0,$$

when $m$ goes to infinity. By applying the first part (DP1) of the dynamic programming principle to $v(t_m, x_m)$ and $\theta_m := \tau_m \wedge (t_m + h_m)$, we get :

$$v(t_m, x_m) \geq E\left[\int_{t_m}^{\theta_m} f(s, X_s^{t_m, x_m}, a)ds + v(\theta_m, X_{\theta_m}^{t_m, x_m})\right].$$

From (3.13), which implies $v \geq v_* \geq \varphi$, we then deduce :

$$\varphi(t_m, x_m) + \gamma_m \geq E\left[\int_{t_m}^{\theta_m} f(s, X_s^{t_m, x_m}, a)ds + \varphi(\theta_m, X_{\theta_m}^{t_m, x_m})\right].$$

We now apply Itô's formula to $\varphi(s, X_s^{t_m, x_m})$ between $t_m$ and $\theta_m$, and we obtain after noting that the stochastic integral term vanishes in expectation due to bounded integrand :

$$\frac{\gamma_m}{h_m} + E\left[\frac{1}{h_m}\int_{t_m}^{\theta_m}\left(-\frac{\partial \varphi}{\partial t} - \mathcal{L}^a \varphi - f\right)(s, X_s^{t_m, x_m}, a)ds\right] \geq 0. \quad (3.14)$$

By the a.s. continuity of the trajectory $X_s^{t_m,x_m}$, we get that for $m$ sufficiently large, $(m \geq N(\omega))$, $\theta_m(\omega) = t_m + h_m$, a.s. We then deduce by the mean-value theorem that the random variable inside the expectation in (3.14) converges a.s. to $-\frac{\partial \varphi}{\partial t}(\bar{t}, \bar{x}) - \mathcal{L}^a \varphi(\bar{t}, \bar{x}) - f(\bar{t}, \bar{x}, a)$ when $m$ goes to infinity. Moreover, since this random variable is bounded by a constant independent of $m$, we may apply the dominated convergence theorem to obtain :

$$-\frac{\partial \varphi}{\partial t}(\bar{t}, \bar{x}) - \mathcal{L}^a \varphi(\bar{t}, \bar{x}) - f(\bar{t}, \bar{x}, a) \geq 0.$$

We conclude from the arbitrariness of $a \in A$. □

**Remark 3.5** The supersolution property of $v$ means that for all $(\bar{t}, \bar{x}) \in [0, T) \times \mathbb{R}^n$ and smooth test function $\varphi$ s.t. $0 = (v_* - \varphi)(\bar{t}, \bar{x}) = \min_{[0,T) \times \mathbb{R}^n}(v_* - \varphi)$, we have

$$-\frac{\partial v}{\partial t}(\bar{t}, \bar{x}) - H(\bar{t}, \bar{x}, D_x \varphi(\bar{t}, \bar{x}), D_x^2 \varphi(\bar{t}, \bar{x})) \geq 0.$$

Recalling condition (3.5), this implies

$$G(\bar{t}, \bar{x}, D_x \varphi(\bar{t}, \bar{x}), D_x^2 \varphi(\bar{t}, \bar{x})) \geq 0,$$

and so

$$\min\left\{-\frac{\partial \varphi}{\partial t}(\bar{t}, \bar{x}) - H(\bar{t}, \bar{x}, D_x \varphi(\bar{t}, \bar{x}), D_x^2 \varphi(\bar{t}, \bar{x})),\right.$$
$$\left. G(\bar{t}, \bar{x}, D_x \varphi(\bar{t}, \bar{x}), D_x^2 \varphi(\bar{t}, \bar{x}))\right\} \geq 0.$$

This is the viscosity supersolution of $v$ to (3.6).

*3.2.2. Viscosity subsolution property*

**Proposition 3.2** *Assume that* (3.5) *holds. Then $v$ is a viscosity subsolution of the HJB variational inequality :*

$$\min\left\{-\frac{\partial v}{\partial t}(t, x) - H(t, x, D_x v(t, x), D_x^2 v(t, x)),\right.$$
$$\left. G(t, x, D_x v(t, x), D_x^2 v(t, x))\right\} = 0, \quad (t, x) \in [0, T) \times \mathbb{R}^n. \quad (3.15)$$

The proof of the subsolution part is based on a contraposition argument and DP2. We then introduce for a given smooth function $\varphi$, the set in $[0, T] \times \mathbb{R}^n$ :

$$\mathcal{M}(\varphi) = \left\{(t, x) \in [0, T] \times \mathbb{R}^n : G(t, x, D_x \varphi(t, x), D_x^2 \varphi(t, x)) > 0\right.$$
$$\left. \text{and} \quad -\frac{\partial \varphi}{\partial t}(t, x) - H(t, x, D_x \varphi(t, x), D_x^2 \varphi(t, x)) > 0\right\}.$$

The following Lemma, which will be also used in the derivation of the terminal condition, is a consequence of DP2 of the dynamic programming principle.

**Lemma 3.1** *Let $\varphi$ be a smooth function on $[0,T] \times \mathbb{R}^n$, and suppose there exist $t_1 < t_2 \leq T$, $\bar{x} \in \mathbb{R}^n$ and $\eta > 0$ s.t. :*

$$[t_1, t_2] \times \bar{B}_\eta(\bar{x}) \in \mathcal{M}(\varphi).$$

*Then,*

$$\sup_{\partial_p([t_1,t_2] \times \bar{B}_\eta(\bar{x}))} (v - \varphi) = \max_{[t_1,t_2] \times \bar{B}_\eta(\bar{x})} (v^* - \varphi),$$

*where $\partial_p([t_1, t_2] \times B_\eta(\bar{x}))$ is the forward parabolic boundary of $[t_1,t_2] \times \bar{B}_\eta(\bar{x})$, i.e. $\partial_p([t_1,t_2] \times \bar{B}_\eta(\bar{x})) = [t_1,t_2] \times \partial \bar{B}_\eta(\bar{x}) \cup \{t_2\} \times \bar{B}_\eta(\bar{x})$.*

Before proving this Lemma, let us show how it immediately implies the required subsolution property.

**Proof of Proposition 3.2.** Let $(\bar{t}, \bar{x}) \in [0,T) \times \mathbb{R}^n$ and $\varphi$ a smooth test function s.t.

$$0 = (v^* - \varphi)(\bar{t}, \bar{x}) = (\text{strict}) \max_{[0,T) \times \mathbb{R}^n} (v^* - \varphi).$$

First, observe that by the continuity condition in (3.5), the set $\mathcal{M}(\varphi)$ is open. Since $(\bar{t}, \bar{x})$ is a strict maximizer of $(v^* - \varphi)$, we then deduce by Lemma 3.1 that $(\bar{t}, \bar{x}) \notin \mathcal{M}(\varphi)$. By definition of $\mathcal{M}(\varphi)$, this means :

$$\min \left\{ -\frac{\partial \varphi}{\partial t}(t,x) - H(\bar{t}, \bar{x}, D_x \varphi(\bar{t}, \bar{x}), D_x^2 \varphi(\bar{t}, \bar{x})), \right.$$
$$\left. G(\bar{t}, \bar{x}, D_x \varphi(\bar{t}, \bar{x}), D_x^2 \varphi(\bar{t}, \bar{x})) \right\} \leq 0,$$

which is the required subsolution inequality.

**Proof of Lemma 3.1.** By definition of $\mathcal{M}(\varphi)$ and $H$, we have for all $a \in A$ :

$$-\frac{\partial \varphi}{\partial t}(t,x) - \mathcal{L}^a \varphi(t,x) - f(t,x,a) > 0, \quad \forall (t,x) \in [t_1,t_2] \times \bar{B}_\eta(\bar{x}) \quad (3.16)$$

We argue by contradiction and suppose on the contrary that :

$$\max_{[t_1,t_2] \times \bar{B}_\eta(\bar{x})} (v^* - \varphi) - \sup_{\partial_p([t_1,t_2] \times \bar{B}_\eta(\bar{x}))} (v - \varphi) := 2\delta.$$

We can choose $(t_0, x_0) \in (t_1, t_2) \times B_\eta(\bar{x})$ s.t. $(v - \varphi)(t_0, x_0) \geq -\delta + \max_{[t_1,t_2] \times \bar{B}_\eta(\bar{x})} (v^* - \varphi)$, and so :

$$(v - \varphi)(t_0, x_0) \geq \delta + \sup_{\partial_p([t_1,t_2] \times \bar{B}_\eta(\bar{x}))} (v - \varphi). \quad (3.17)$$

Fix now $\varepsilon = \delta/2$, and apply DP2 to $v(t_0, x_0)$ : there exists $\hat{\alpha}^\varepsilon \in \mathcal{A}$ s.t.

$$v(t_0, x_0) - \varepsilon \leq E\left[ \int_{t_0}^{\theta} f(s, X_s^{t_0,x_0}, \hat{\alpha}_s^\varepsilon) ds + v(\theta, X_\theta^{t_0,x_0}) \right], \quad (3.18)$$

where we choose

$$\theta = \inf\left\{s \geq t_0 : (s, X_s^{t_0,x_0}) \notin [t_1, t_2] \times \bar{B}_\eta(\bar{x})\right\}.$$

First, notice that by continuity of $X^{t_0,x_0}$, we have $(\theta, X_\theta^{t_0,x_0}) \in \partial_p([t_1, t_2] \times B_\eta(\bar{x}))$. Since from (3.17), we have $v \leq \varphi + (v - \varphi)(t_0, x_0) - \delta$ on $\partial_p([t_1, t_2] \times B_\eta(\bar{x}))$, we get with (3.18) :

$$-\varepsilon \leq E\left[\int_{t_0}^\theta f(s, X_s^{t_0,x_0}, \hat{\alpha}_s^\varepsilon)ds + \varphi(\theta, X_\theta^{t_0,x_0}) - \varphi(t_0, x_0)\right] - \delta.$$

Applying Itô's formula to $\varphi(s, X_s^{t_0,x_0})$ between $s = t_0$ and $s = \theta$, we obtain :

$$E\left[\int_{t_0}^\theta \left(-\frac{\partial\varphi}{\partial t}(s, X_s^{t_0,x_0}) - \mathcal{L}^{\hat{\alpha}_s^\varepsilon}\varphi(s, X_s^{t_0,x_0}) - f(s, X_s^{t_0,x_0}, \hat{\alpha}_s^\varepsilon)\right)ds\right] \leq \varepsilon - \delta.$$

Since, by definition of $\theta$, $(s, X_s^{t_0,x_0})$ lies in $[t_1, t_2] \times \bar{B}_\eta(\bar{x})$ for all $t_0 \leq s \leq \theta$, we get with (3.16) the required contradiction : $0 \leq \varepsilon - \delta = -\delta/2$.

*3.2.3. Terminal condition*

We start with the following Lemma.

**Lemma 3.2** *Suppose that $f$ and $g$ are lower-bounded or satisfy a linear growth condition, and $g$ is lowersemicontinuous. Then,*

$$\underline{v}(x) \geq g(x), \quad \forall x \in \mathbb{R}^n.$$

**Proof.** Take some arbitrary sequence $(t_m, x_m) \to (T, x)$ with $t_m < T$ and fix some constant control $\alpha \equiv a \in \mathcal{A}$. By definition of the value function, we have :

$$v(t_m, x_m) \geq E\left[\int_{t_m}^T f(s, X_s^{t_m,x_m}, a)ds + g(X_T^{t_m,x_m})\right].$$

Under the linear growth or lower-boundeness condition on $f$ and $g$, we may apply the dominated convergence theorem or Fatou's lemma, and so :

$$\liminf_{m\to\infty} v(t_m, x_m) \geq E\left[\liminf_{m\to\infty} g(X_T^{t_m,x_m})\right]$$
$$\geq g(x),$$

by the lowersemicontinuity of $g$ and the continuity of the flow $X_T^{t,x}$ in $(t, x)$. □

The supersolution property (3.8) for the terminal condition is then obtained with the following result.

**Lemma 3.3** *Under (3.5), $\underline{v}$ is a viscosity supersolution of :*

$$G(T, x, D_x\underline{v}(x), D_x^2\underline{v}(x)) = 0, \quad on \ \mathbb{R}^n.$$

**Proof.** Let $\bar{x} \in \mathbb{R}^n$ and $\psi$ a smooth function on $\mathbb{R}^n$ s.t.

$$0 = (\underline{v} - \psi)(\bar{x}) = \min_{\mathbb{R}^n}(\underline{v} - \psi). \qquad (3.19)$$

By definition of $\underline{v}$, there exists a sequence $(s_m, y_m)$ converging to $(T, \bar{x})$ with $s_m < T$ and

$$\lim_{m \to \infty} v_*(s_m, y_m) = \underline{v}(\bar{x}). \qquad (3.20)$$

Consider the auxiliary test function:

$$\varphi_m(t, x) = \psi(x) - |x - \bar{x}|^4 + \frac{T - t}{(T - s_m)^2},$$

and choose $(t_m, x_m) \in [s_m, T] \times \bar{B}_1(\bar{x})$ as a minimum of $(v_* - \varphi_m)$ on $[s_m, T] \times \bar{B}_1(\bar{x})$.

*Step 1.* We claim that, for sufficiently large $m$, $t_m < T$ and $x_m$ converges to $\bar{x}$, so that $(t_m, x_m)$ is a local minimizer of $(v_* - \varphi_m)$. Indeed, recalling $\underline{v}(\bar{x}) = \psi(\bar{x})$ and (3.20), we have for sufficiently large $m$:

$$(v_* - \varphi_m)(s_m, y_m) \leq -\frac{1}{2(T - s_m)} < 0. \qquad (3.21)$$

On the other hand, for any $x \in \mathbb{R}^n$, we have

$$(v_* - \varphi_m)(T, x) = \underline{v}(x) - \psi(x) + |x - \bar{x}|^4 \geq \underline{v}(x) - \psi(x) \geq 0, \quad (3.22)$$

by (3.19). The two inequalities (3.21)-(3.22) show that $t_m < T$ for large $m$. We can suppose that $x_m$ converges, up to a subsequence, to some $x_0 \in \bar{B}_1(\bar{x})$. From (3.19), since $s_m \leq t_m$ and $(t_m, x_m)$ is a minimum of $(v_* - \psi_m)$, we have:

$$\begin{aligned}
0 &\leq (\underline{v} - \psi)(x_0) - (\underline{v} - \psi)(\bar{x}) \\
&\leq \liminf_{m \to \infty} \left[(v_* - \varphi_m)(t_m, x_m) - (v_* - \varphi_m)(s_m, y_m) - |x_m - \bar{x}|^4\right] \\
&\leq -|x_0 - \bar{x}|^4,
\end{aligned}$$

which proves that $x_0 = \bar{x}$.

*Step 2.* Since $(t_m, x_m)$ is a local minimizer of $(v_* - \varphi_m)$, the viscosity supersolution property of $v_*$ holds at $(t_m, x_m)$ with the test function $\varphi_m$, and so for every $m$:

$$G(t_m, x_m, D_x\varphi_m(t_m, x_m), D_x^2\varphi_m(t_m, x_m)) \geq 0. \qquad (3.23)$$

Now, since $D_x\varphi_m(t_m, x_m) = D\psi(x_m) - 4(x_m - \bar{x})|x_m - \bar{x}|^2$, $D_x^2\varphi_m(t_m, x_m) = D^2\psi(x_m) - 4|x_m - \bar{x}|^2 I_n - 4(x_m - x)(x_m - \bar{x})'$, recalling that $G$ is continuous, and $(t_m, x_m)$ converges to $(T, \bar{x})$, we get from (3.23):

$$G(T, \bar{x}, D\psi(\bar{x}), D^2\psi(\bar{x})) \geq 0.$$

This is the required supersolution inequality. □

We finally turn to the subsolution property for the terminal condition.

**Lemma 3.4** *Suppose that $g$ is uppersemicontinuous and (3.5) holds. Then, $\bar{v}$ is a viscosity subsolution of :*

$$\min\left[\bar{v}(x) - g(x)\,,\, G(T, x, D\bar{v}(x), D^2\bar{v}(x))\right] \;=\; 0, \quad on\; \mathbb{R}^n.$$

**Proof.** Let $\bar{x} \in \mathbb{R}^n$ and $\psi$ a smooth function on $\mathbb{R}^n$ s.t.

$$0 \;=\; (\bar{v} - \psi)(\bar{x}) \;=\; \max_{\mathbb{R}^n}(\bar{v} - \psi). \tag{3.24}$$

We have to show that whenever

$$\bar{v}(\bar{x}) \;>\; g(\bar{x}), \tag{3.25}$$

then

$$G(T, \bar{x}, D\psi(\bar{x}), D^2\psi(\bar{x})) \;\leq\; 0. \tag{3.26}$$

So, suppose that (3.25) holds, and let us consider the auxiliary test function :

$$\varphi_m(t, x) \;=\; \psi(x) + |x - \bar{x}|^4 + m(T - t).$$

We argue by contradiction and suppose on the contrary that

$$G(T, \bar{x}, D\psi(\bar{x}), D^2\psi(\bar{x})) \;>\; 0.$$

Since $D_x\varphi_m(t, x) = D\psi(x) - 4(x - \bar{x})|x - \bar{x}|^2 \to D\psi(\bar{x})$, $D_x^2\varphi_m(t, x) = D^2\psi(x) - 4I_n|x - \bar{x}|^2 - 4(x - \bar{x})(x - \bar{x})' \to D^2\psi(\bar{x})$ when $x$ tends to $\bar{x}$, there exists, by continuity of $G$, $s_0 < T$ and $\eta > 0$ s.t. for all $m$ :

$$G(t, x, D_x\varphi_m(t, x), D_x^2\varphi_m(t, x)) \;>\; 0, \quad \forall (t, x) \in [s_0, T] \times \bar{B}_\eta(\bar{x}). \tag{3.27}$$

Under condition (3.5), the function $H(t, x, D_x\varphi_m(t, x), D_x^2\varphi_m(t, x))$ is then finite on the compact set $[s_0, T] \times \bar{B}_\eta(\bar{x})$ and by continuity of $H$ on $\text{int}(\text{dom}(H))$, there exists some constant $h_0$ (independent of $m$) s.t.

$$H(t, x, D_x\varphi_m(t, x), D_x^2\varphi_m(t, x)) \;\leq\; h_0, \quad \forall (t, x) \in [s_0, T] \times \bar{B}_\eta(\bar{x}). \tag{3.28}$$

<u>Step 1.</u> Since by definition $\bar{v} \geq \underline{v}$, we have from Lemma 3.2 :

$$\bar{v} \;\geq\; g. \tag{3.29}$$

Hence, for all $x \in \mathbb{R}^n$,

$$(v - \varphi_m)(T, x) \;=\; (g - \psi)(x) - |x - \bar{x}|^4 \;\leq\; (\bar{v} - \psi)(x) - |x - \bar{x}|^4$$
$$\leq\; -|x - \bar{x}|^4 \;\leq\; 0 \tag{3.30}$$

by (3.24). This implies : $\sup_{B_\eta(\bar{x})}(v - \varphi_m)(T, .) \leq 0$. We claim that

$$\limsup_{m \to \infty} \sup_{B_\eta(\bar{x})} (v - \varphi_m)(T, .) \;<\; 0. \tag{3.31}$$

On the contrary, there exists a subsequence of $(\varphi_m)$, still denoted $(\varphi_m)$ s.t. :

$$\lim_{m\to\infty} \sup_{B_\eta(\bar{x})} (v - \varphi_m)(T, .) = 0,$$

For each $m$, let $(x_m^k)_k$ be a maximizing sequence of $(v - \varphi_m)(T, .)$ on $B_\eta(\bar{x})$, i.e.

$$\lim_{m\to\infty} \lim_{k\to\infty} (v - \varphi_m)(T, x_m^k) = 0.$$

Now, from (3.30), we have $(v - \varphi_m)(T, x_m^k) \leq -|x_m^k - \bar{x}|^4$, which combined with the above equality shows that :

$$\lim_{m\to\infty} \lim_{k\to\infty} x_m^k = \bar{x}.$$

Hence,

$$\begin{aligned} 0 &= \lim_{m\to\infty} \lim_{k\to\infty} (v - \varphi_m)(T, x_m^k) = \lim_{m\to\infty} \lim_{k\to\infty} g(x_m^k) - \psi(\bar{x}) \\ &\leq g(\bar{x}) - \psi(\bar{x}) < (\bar{v} - \psi)(\bar{x}), \end{aligned}$$

by the uppersemicontinuty of $g$ and (3.25). This contradicts $(\bar{v} - \psi)(\bar{x}) = 0$ in (3.24).

*Step 2.* Take a sequence $(s_m)$ converging to $T$ with $s_0 \leq s_m < T$. Let us consider a maximizing sequence $(t_m, x_m)$ of $v^* - \varphi_m$ on $[s_m, T] \times \partial \bar{B}_\eta(\bar{x})$. Then

$$\limsup_{m\to\infty} \sup_{[s_m,T]\times\partial\bar{B}_\eta(\bar{x})} (v^* - \varphi_m) \leq \limsup_{m\to\infty} (v^*(t_m, x_m) - \psi(x_m)) - \eta^4.$$

Since $t_m$ converges to $T$ and $x_m$, up to a subsequence, converges to some $x_0 \in \partial \bar{B}_\eta(\bar{x})$, we have by definition of $\bar{v}$ :

$$\limsup_{m\to\infty} \sup_{[s_m,T]\times\partial\bar{B}_\eta(\bar{x})} (v^* - \varphi_m) \leq (\bar{v} - \psi)(x_0) - \eta^4 \leq -\eta^4, \quad (3.32)$$

by (3.24). Recall also from (3.24) that $(v^* - \varphi_m)(T, \bar{x}) = (\bar{v} - \psi)(\bar{x}) = 0$. Therefore, with (3.31) and (3.32), we deduce that for $m$ large enough :

$$\sup_{[s_m,T]\times\partial\bar{B}_\eta(\bar{x})} (v - \psi_m) < 0 = (v^* - \varphi_m)(T, \bar{x}) \leq \max_{[s_m,T]\times\partial B_\eta(\bar{x})} (v^* - \varphi_m).$$

In view of Lemma 3.1, this proves that for $m$ large enough,

$$[s_m, T] \times \bar{B}_\eta(\bar{x}) \quad \text{is not a subset of} \quad \mathcal{M}(\varphi_m). \tag{3.33}$$

*Step 3.* From (3.28), notice that for all $(t, x) \in [s_m, T] \times \bar{B}_\eta(\bar{x})$, we have :

$$-\frac{\partial \varphi_m}{\partial t}(t, x) - H(t, x, D_x\varphi_m(t, x), D_x^2\varphi_m(t, x)) \geq m - h_0 > 0$$

for $m$ large enough. In view of (3.33) and by definition of $\mathcal{M}(\varphi_m)$, we then may find some element $(t, x) \in [s_m, T] \times \bar{B}_\eta(\bar{x})$ s.t.

$$G(t, x, D_x\varphi_m(t, x), D_x^2\varphi_m(t, x)) \leq 0.$$

This is in contradiction with (3.27). □

### 3.3. Some applications in finance

*3.3.1. Smooth-fit property of one-dimensional convex singular problem*

We consider the singular control problem introduced in Remark 2.2 with linear dynamics :

$$dX_t = \beta X_t dt + \gamma X_t dt + \alpha_t dt, \qquad (3.34)$$

where $\beta$ and $\gamma$ are constants with $\gamma > 0$, and $\alpha$ is the control valued in $\mathbb{R}_+$. We consider the infinite horizon problem :

$$v(x) = \sup_{\alpha \in \mathcal{A}} E\left[\int_0^\infty e^{-\rho t}\left(\hat{f}(X_t^x) - \lambda \alpha_t\right) dt\right],$$

where $X_t^x$ is the solution to (3.34) starting from $x$ at time 0. It is convenient to redefine control in terms $L_t = \int_0^t \alpha_s ds$, and in fact to enlarge the set of control processes to the set of nondecreasing, cadlag adapted processes $L$. This allows displacement in time of the control and also ensures the existence of an optimal control. This singular control problem models the irreversible investment problem for a firm. The process $X$ represents the production capacity of a firm, which may increase it by transferring capital from another sector activity. $\beta \leq 0$ is the depreciation rate of the capacity production, $\gamma$ its volatility, $L_t$ is the cumulated number of capital received by the firm until time $t$ for a cost $\int_0^t \lambda dL_t$, with $\lambda > 0$, interpreted as a conversion factor from an activity sector to another one. $\hat{f}$ is the running profit function of the firm, assumed to be concave and finite, i.e. of linear growth condition, on $(0, \infty)$.

The Hamiltonian of this one-dimensional convex singular control problem is :

$$H(x, p, M) = \begin{cases} \beta x p + \frac{1}{2}\gamma^2 x^2 M + \hat{f}(x) & \text{if } \lambda - p \geq 0 \\ \infty & \text{if } \lambda - p < 0. \end{cases}$$

Hence, the associated HJB variational inequality takes the form :

$$\min\left[\rho v(x) - \beta x v'(x) - \frac{1}{2}\gamma^2 x^2 v''(x) - \hat{f}(x) \, , \, \lambda - v'(x)\right] = 0. \quad (3.35)$$

One-dimensional convex singular control problems were largely sudied in the 1980's and in particular their connection with optimal stopping problems. Formally, the derivative of the value function for a singular control problem is the value function of an optimal stopping problem. In view of the smooth-fit principle for optimal stopping problems, which states that the associated value function should be differentiable, it is expected that the value function for a singular control problem should be twice differentiable. We show here how to derive simply the twice continuous differentiability for the value function of the above one-dimensional convex singular problem by using a viscosity solutions argument.

First, we easily check that by the linearity of the dynamics of the controlled process $X$ and the concavity and linear growth condition of $\hat{f}$, the value function $v$ is concave with linear growth condition on $(0, \infty)$. Hence, $v$ is continuous on $(0, \infty)$, and admits a left and right derivative $v'_-(x)$ and $v'_+(x)$ for each $x > 0$, with $v'_+(x) \leq v'_-(x)$. We also see that

$$v'_-(x) \leq \lambda, \quad x > 0. \tag{3.36}$$

Indeed, fix $x > 0$, and let us consider for any arbitrary control $L$ and positive number $0 < l < x$, the control $\tilde{L}$ consisting of an immediate jump of $L$ of size $l$ at time $0$. Then, the associated state process $\tilde{X}$ satisfies: $\tilde{X}^{x-l} = X^x$ and by definition of the value function, we have:

$$\begin{aligned} v(x-l) &\geq E\left[\int_0^\infty e^{-\rho t}\left(\hat{f}(\tilde{X}_t^{x-l})dt - \lambda d\tilde{L}_t\right)\right] \\ &= E\left[\int_0^\infty e^{-\rho t}\left(\hat{f}(X_t^x)dt - \lambda dL_t\right)\right] - \lambda l. \end{aligned}$$

From the arbitrariness of $L$, this implies $v(x-l) \geq v(x) - \lambda l$ and so the required result (3.36) by sending $l$ to zero. By setting $x_b = \inf\{x \geq 0 : v'_+(x) < \lambda\}$, and using the concavity of $v$, we also easily see that

$$\mathcal{NT} := \{x > 0 : v'_-(x) < \lambda\} = (x_b, \infty), \tag{3.37}$$

$v$ is differentiable on $(0, x_b)$ with

$$v'(x) = \lambda, \quad \text{on } (0, x_b). \tag{3.38}$$

We can also easily check that $0 < x_b < \infty$.

According to Theorem 3.1, we know that $v$ is viscosity solution to the HJB variational inequality (3.35). By exploiting furthermore the concavity of $v$ in this one-dimensional problem, we show that $v$ is actually a classical smooth solution.

**Theorem 3.3** *The value function $v$ is a classical solution in $C^2((0, \infty))$ of (3.35).*

**Proof.** *Step 1.* We first prove that $v$ in $C^1$ on $(0, \infty)$. We argue by contradiction and suppose then on the contrary that $v'_+(\bar{x}) < v'_-(\bar{x})$ for some $\bar{x} > 0$. Let us fix $\bar{p} \in (v'_+(\bar{x}), v'_-(\bar{x}))$ and consider the smooth function

$$\varphi_\varepsilon(x) = v(\bar{x}) + \bar{p}(x - \bar{x}) - \frac{1}{2\varepsilon}(x - \bar{x})^2,$$

with $\varepsilon > 0$. Observe that $\bar{x}$ is a local maximum of $(v - \varphi_\varepsilon)$ with $\varphi_\varepsilon(\bar{x}) = v(\bar{x})$. Since $\varphi'_\varepsilon(\bar{x}) = \bar{p} < \lambda$ by (3.36) and $\varphi''_\varepsilon(\bar{x}) = 1/\varepsilon$, the viscosity subsolution property implies:

$$\rho\varphi(\bar{x}) - \beta\bar{x}\bar{p} + \frac{1}{2\varepsilon}\gamma^2\bar{x}^2 - \hat{f}(\bar{x}) \leq 0. \tag{3.39}$$

For $\varepsilon$ sufficiently small, we get the required contradiction, and so $v'_+(\bar{x}) = v'_-(\bar{x})$.

*Step 2.* By (3.38), $v$ is $C^2$ on $(0, x_b)$ and satisfies $v'(x) = \lambda$, $x \in (0, x_b)$. From Step 1, we have $\mathcal{NT} = (x_b, \infty) = \{x > 0 : v'(x) < \lambda\}$. Let us check that $v$ is a viscosity solution of :

$$\rho v(x) - \beta x v'(x) - \frac{1}{2}\gamma^2 x^2 v''(x) - \hat{f}(x) = 0, \quad x \in (x_b, \infty). \quad (3.40)$$

Let $\bar{x} \in (x_b, \infty)$ and $\varphi$ a smooth test function $C^2$ on $(x_b, \infty)$ s.t. $\bar{x}$ is a local maximum of $v - \varphi$, with $(v - \varphi)(\bar{x}) = 0$. Since $\varphi'(\bar{x}) = v'(\bar{x}) < \lambda$, the viscosity subsolution property of $v$ to (3.35) implies :

$$\rho \varphi(\bar{x}) - \beta \bar{x} \varphi'(\bar{x}) - \frac{1}{2}\gamma^2 \bar{x}^2 \varphi''(\bar{x}) - \hat{f}(\bar{x}) \leq 0.$$

This shows that $v$ is a viscosity subsolution of (3.40) on $(x_b, \infty)$. The proof of supersolution viscosity of (3.40) is similar. Consider now arbitrary $x_1 < x_2$ in $(x_b, \infty)$, and the Dirichlet problem :

$$\rho V(x) - \beta x V'(x) - \frac{1}{2}\gamma^2 x^2 V''(x) - \hat{f}(x) = 0, \quad x \in (x_1, x_2) \quad (3.41)$$

$$V(x_1) = v(x_1), \qquad V(x_2) = v(x_2). \quad (3.42)$$

Classical results provide the existence and uniqueness of a smooth solution $V$, $C^2$ on $(x_1, x_2)$, to (3.41)-(3.42). In particular, this smooth solution $V$ is a viscosity solution to (3.40) on $(x_1, x_2)$. By standard comparison principle for viscosity solutions (here for linear PDE in bounded domain), we deduce that $v = V$ on $(x_1, x_2)$. From the arbitrariness of $x_1$ and $x_2$, this proves that $v$ is $C^2$ on $(x_b, \infty)$ and satisfies (3.40) in the classical sense.

*Step 3.* Recall that $v$ is $C^2$ on $(0, x_b)$ with $v'(x) = \lambda$. It then remains to prove the $C^2$ condition of $v$ in $x_b$. Let $x \in (0, x_b)$. The viscosity subsolution property (3.35) of $v$ in $x$ applied with the smooth test function $\varphi = v$ on $(0, x_b)$, where $\varphi'(x) = \lambda$, $\varphi''(x) = 0$, implies that $v$ satisfies

$$\rho v(x) - \beta \lambda x - \hat{f}(x) \geq 0, \quad 0 < x < x_b,$$

in the classical sense. By sending $x$ to $x_b$, we get :

$$\rho v(x_b) - \beta \lambda x_b - \hat{f}(x_b) \geq 0. \quad (3.43)$$

On the other hand, by the $C^1$ condition of $v$ in $x_b$, we have by sending $x$ to $x_b$ in (3.40) :

$$\rho v(x_b) - \beta \lambda x_b - \hat{f}(x_b) = \frac{1}{2}\gamma^2 x_b^2 v''(x_b^+). \quad (3.44)$$

From the concavity of $v$, the r.h.s. of (3.44) is nonpositive, which combined with (3.43), shows that $v''(x_b^+) = 0$. This proves that $v$ is $C^2$ in $x_b$ with

$$v''(x_b) = 0.$$

□

Once we have the $C^2$ regularity of the value function, we complete its explicit characterization as follows. We solve the ode (3.38) on $(0, x_b)$ and the ode (3.40) in $(x_b, \infty)$. We have then four unknown parameters : one coming from the first order ode (3.38), two from the second order ode (3.40), and the parameter $x_b$. These four parameters are explicitly determined by the linear growth condition of $v$ and the $C^2$ condition of $v$ in $x_b$. Moreover, the optimal state process and control are given by the reflected process $\hat{X}$ at the lower level $x_b$ and by the local time $\hat{L}$ at $x_b$. We refer to [46] for the details in a more general framework.

The concavity property of $v$, resulting from the linear dynamics of $X$, is crucial for the $C^2$ regularity of the value function : there are examples where the value function is not convex and not twice continuously differentiable. Some extensions to the two-dimensional case are studied in [84] for a singular control problem arising from a transaction costs model, and the value function is shown to be twice differentiable under some assumptions. Singular control problems in multidimensional transaction costs model were recently studied in [50] but no regularity results are stated : only a viscosity characterization is given. It is an open question to know whether regularity results can be extended in such a multidimensional context.

*3.3.2. Superreplication in uncertain volatility model*

We consider the controlled diffusion :

$$dX_s = \alpha_s X_s dW_s,$$

valued in $(0, \infty)$ (for an initial condition $x > 0$), and where the control $\alpha$ is valued in $A = [\underline{a}, \bar{a}]$, $0 \leq \underline{a} \leq \bar{a} \leq \infty$. To avoid trivial cases, we assume $\bar{a} > 0$ and $\underline{a} \neq \infty$. In finance, $\alpha$ represents the uncertain volatility process of the asset price $X$. Given a continuous function $g$ with linear growth condition on $(0, \infty)$, representing the payoff function of an European option, we are interested in the calculation of its superreplication price given by, see e.g. [34] :

$$v(t,x) = \sup_{\alpha \in \mathcal{A}} E\left[g(X_T^{t,x})\right], \quad (t,x) \in [0,T] \times (0, \infty) \qquad (3.45)$$

Since the process $X$ is a nonnegative supermartingale, it is easy to see that $v$ inherits from $g$ the linear growth condition. In particular, $v$ is locally bounded.

The Hamiltonian of this stochastic control problem is

$$H(x, M) = \sup_{a \in [\underline{a}, \bar{a}]} \left\{ \frac{1}{2} a^2 x^2 M \right\}, \quad (x, M) \in (0, \infty) \times \mathbb{R}.$$

We shall then distinguish two cases according to the finiteness of the volatility upper bound $\bar{a}$.

**Bounded volatility**

We suppose

$$\bar{a} < \infty.$$

In this regular case, the Hamiltonian $H$ is finite on the whole domain $(0, \infty) \times \mathbb{R}$, and is explicitly given by:

$$H(x, M) = \frac{1}{2}\hat{a}^2(M)x^2 M,$$

with

$$\hat{a}(M) = \begin{cases} \bar{a} & \text{if } M \geq 0 \\ \underline{a} & \text{if } M < 0. \end{cases}$$

According to the results in paragraph 3.2 (Theorems 3.1, 3.2, and their following remarks), we have then the characterization result on the superreplication price $v$:

**Theorem 3.4** *Suppose $\bar{a} < \infty$. Then $v$ is continuous on $[0, T] \times (0, \infty)$ and is the unique viscosity solution with linear growth condition of the so-called Black-Scholes-Barenblatt equation:*

$$-\frac{\partial v}{\partial t} - \frac{1}{2}\hat{a}^2\left(\frac{\partial^2 v}{\partial x^2}\right)x^2\frac{\partial^2 v}{\partial x^2} = 0, \quad (t, x) \in [0, T) \times (0, \infty), \quad (3.46)$$

*satisfying the terminal condition:*

$$v(T, x) = g(x), \quad x \in (0, \infty). \quad (3.47)$$

**Remark 3.6** When $\underline{a} > 0$, there is existence and uniqueness of a smooth solution to the Black-Scholes-Barenblatt PDE (3.46) together with the Cauchy condition (3.47), and so $v$ is this smooth solution.

**Unbounded volatility**

We suppose

$$\bar{a} = \infty.$$

In this singular case, the Hamiltonian is given by:

$$H(x, M) = \begin{cases} \frac{1}{2}\underline{a}^2 x^2 M & \text{if } -M \geq 0 \\ \infty & \text{if } -M < 0. \end{cases}$$

According to Theorem 3.1, $v$ is then a viscosity solution of:

$$\min\left\{-\frac{\partial v}{\partial t} - \frac{1}{2}\underline{a}^2 x^2 \frac{\partial^2 v}{\partial x^2}, -\frac{\partial^2 v}{\partial x^2}\right\} = 0, \quad (t, x) \in [0, T) \times (0, \infty). (3.48)$$

Moreover, from Theorem 3.2 and its following remark, the terminal condition associated to the variational inequality (3.48) is given by:

$$v(T^-, x) = \hat{g}(x), \quad x \in (0, \infty),$$

where $\hat{g}$ is the smallest function above $g$ and satisfying in the viscosity sense

$$-\frac{\partial^2 \hat{g}}{\partial x^2}(x) \geq 0 \quad \text{on } (0, \infty). \tag{3.49}$$

If $\hat{g}$ were smooth, the previous inequality (3.49) characterizes the concavity of $\hat{g}$. Actually, this is still true when (3.49) holds only in the viscosity sense, see [24]. Therefore, $\hat{g}$ is the upper concave envelope of $g$. We can then explicitly characterize the superreplication price $v$.

**Theorem 3.5** *Suppose $\bar{a} = \infty$. Then $v = w$ on $[0, T) \times (0, \infty)$ where $w$ is the Black-Scholes price for the payoff function $\hat{g}(x)$ :*

$$w(t, x) = E\left[\hat{g}\left(\hat{X}_T^{t,x}\right)\right], \quad \forall (t, x) \in [0, T] \times (0, \infty), \tag{3.50}$$

*in a Black-Scholes model with lower volatility $\underline{a}$, i.e. $\{\hat{X}_s^{t,x}, t \leq s \leq T\}$ is the solution to*

$$d\hat{X}_s = \underline{a}\hat{X}_s dW_s, \quad t \leq s \leq T, \quad \hat{X}_t = x.$$

**Proof.** First, observe that $w$ is the solution to the Black-Scholes PDE :

$$-\frac{\partial w}{\partial t} - \frac{1}{2}\underline{a}^2 x^2 \frac{\partial^2 w}{\partial x^2} = 0, \quad \text{on } [0, T) \times (0, \infty), \tag{3.51}$$

with the terminal condition :

$$w(T^-, x) = w(T, x) = \hat{g}(x).$$

Since $\hat{g}$ is concave and it is well-known that Black-Scholes price inherits concavity from its payoff, we deduce that

$$-\frac{\partial^2 w}{\partial x^2}(t, x) \geq 0 \quad \text{on } [0, T) \times (0, \infty).$$

Together with the previous equality (3.51), this proves that $w$ is solution to :

$$\min\left\{-\frac{\partial w}{\partial t} - \frac{1}{2}\underline{a}^2 x^2 \frac{\partial^2 w}{\partial x^2}, -\frac{\partial^2 w}{\partial x^2}\right\} = 0, \quad (t, x) \in [0, T) \times (0, \infty).$$

Finally, since $w$ and $v$ satisfy the same terminal condition, we conclude with a comparison principle on the PDE (3.48). □

Other variations and extensions to financial models with portfolio constraints or transaction costs are studied in [24] and [25] for the calculation of superreplication price by a viscosity solutions approach.

## 4. From stochastic control to backward stochastic differential equations

As already mentioned in paragraph 2.3, there is a strong relation between the stochastic maximum principle and backward stochastic differential equations, BSDE in short. A BSDE is usually written in the form

$$dY_t = -f(t, Y_t, Z_t)dt + Z_t dW_t, \quad Y_T = \xi. \tag{4.1}$$

Here $\mathbb{F} = (\mathcal{F}_t)_{0 \leq t \leq T}$ is the filtration generated by $W$, the terminal condition $\zeta$ is an $\mathcal{F}_T$-measurable random variable, and a solution to (4.1) is a pair of $\mathbb{F}$-adapted processes $(Y, Z)$ s.t.

$$Y_t = \zeta + \int_t^T f(s, Y_s, Z_s)ds - \int_t^T Z_s dW_s, \quad 0 \leq t \leq T.$$

Existence and uniqueness of a solution $(Y, Z)$ are usually stated under a Lipschitz condition on $f$, see Pardoux and Peng [73]. Some extensions to the quadratic growth case are studied in [52] and more recently in [42]. We refer to [32] for an account on the theory of BSDE and to [33] for its applications in finance.

### 4.1. Stochastic maximum principle and BSDE

We state a sufficient maximum principle for a stochastic control problem in the framework described in Section 2 : a controlled state process

$$dX_s = b(X_s, \alpha_s)ds + \sigma(X_s, \alpha_s)dW_s, \quad X_0 = x,$$

and a gain functional to maximize :

$$J(\alpha) = E\left[\int_0^T f(t, X_t, \alpha_t)dt + g(X_T)\right].$$

Moreover, the coefficients $b(t, x, \omega)$, $\sigma(t, x, \omega)$ of the state process as well as the gain functions $f(t, x, \omega)$ and $g(x, \omega)$ may depend on $\omega$ : $b(t, x, \omega)$, $\sigma(t, x, \omega)$ and $f(t, x, \omega)$ are $\mathcal{F}_t$-measurable for all $(t, x)$ and $g(x, \omega)$ is $\mathcal{F}_T$-measurable. For simplicity of notation, we omit the dependence in $\omega$. We consider the generalized Hamiltonian $\mathcal{H} : [0, T] \times \mathbb{R}^n \times A \times \mathbb{R}^n \times \mathbb{R}^{n \times d} \times \Omega \to \mathbb{R}$ :

$$\mathcal{H}(t, x, a, y, z) = b(x, a).y + \text{tr}(\sigma'(x, a)z) + f(t, x, a). \tag{4.2}$$

We also omit the dependence of $\mathcal{H}$ in $\omega$. Notice the difference with the Hamiltonian (2.7) introduced in the dynamic programming approach. We suppose that $\mathcal{H}$ is differentiable in $x$ with derivative denoted $D_x\mathcal{H}$, and for each $\alpha \in \mathcal{A}$, we consider the BSDE, also called the adjoint equation :

$$dY_t = -D_x\mathcal{H}(t, X_t, \alpha_t, Y_t, Z_t)dt + Z_t dW_t, \quad Y_T = D_x g(X_T). \tag{4.3}$$

The following sufficiency theorem is proved e.g. in [89].

**Theorem 4.1** *Assume that $g$ is concave in $x$. Let $\hat{\alpha} \in \mathcal{A}$ and $\hat{X}$ the associated controlled state process. Suppose there exists a solution $(\hat{Y}, \hat{Z})$ to the corresponding BSDE* (4.3) *such that almost surely:*

$$\mathcal{H}(t, \hat{X}_t, \hat{\alpha}_t, \hat{Y}_t, \hat{Z}_t) = \max_{a \in A} \mathcal{H}(t, \hat{X}_t, a, \hat{Y}_t, \hat{Z}_t), \quad 0 \leq t \leq T, \qquad (4.4)$$

*and*

$$(x, a) \rightarrow \mathcal{H}(t, x, a, \hat{Y}_t, \hat{Z}_t) \quad \text{is a concave function}, \qquad (4.5)$$

*for all $t \in [0, T]$. Then $\hat{\alpha}$ is an optimal control, i.e.*

$$J(\hat{\alpha}) = \sup_{\alpha \in \mathcal{A}} J(\alpha).$$

## *4.2. Application : linear-quadratic control problem and mean-variance hedging*

Mean-variance hedging is a popular criterion in finance : this is the problem of finding for a given option payoff the best approximation by means of a self-financed wealth processes, and the optimality criterion is the expected square error. This problem has been solved in high generality for continuous semimartingale price process and general filtration by martingales and projection techniques, see [76] and [83] for recent overviews. On the other hand, in a series of recent papers, this problem has been reformulated and treated as a linear-quadratic (LQ) control problem at increasing levels of generality, mostly for Itô processes and Brownian filtration, see e.g. [54], [53], [67], [14]. When the coefficients are random, which generalize the classical case of deterministic linear-quadratic control problem, the adjoint equations lead to a system of BSDEs. This gives some new insight and results on general LQ control problem, and provides also fairly explicit results for the mean-variance hedging.

The problem is formulated as follows. For simplicity of notation, we consider the one-dimensional case. The dynamics of the linear controlled state process is :

$$dX_t = (r_t X_t + \alpha_t \mu_t)dt + \sigma_t \alpha_t dW_t, \quad X_0 = x,$$

and the objective is to minimize over controls $\alpha$ valued in $\mathbb{R}$, for a given square integrable $\mathcal{F}_T$-measurable random variable $\xi$, the quadratic functional :

$$J(\alpha) = E[X_T - \xi]^2.$$

The coefficients $r_t$, $\mu_t$ and $\sigma_t$ are assumed to be bounded and adapted with respect to $\mathbb{F}$, the natural filtration of $W$. We also assume a uniform ellipticity condition on $\sigma$ : $\sigma_t \geq \varepsilon$ a.s. for some $\varepsilon > 0$. In finance, $X$ is the wealth process of a self-financed portfolio $\alpha$, representing the amount invested in a stock asset of excess return $\mu_t$ with respect to the interest rate $r_t$, and volatility $\sigma_t$. $\xi$ is the payoff at time $T$ of an option that we are trying to approximate by $X$ according to the quadratic error criterion.

We show how to solve this optimization problem by applying the stochastic minimum principle in Theorem 4.1. Observe that the generalized Hamiltonian for this stochastic LQ problem takes the form

$$\mathcal{H}(x,a,y,z) = (r_t x + a\mu_t)y + \sigma_t a z,$$

and the adjoint BSDE (4.3) is written, given $\alpha \in \mathcal{A}$ :

$$dY_t = -r_t Y_t dt + Z_t dW_t, \quad Y_T = 2(X_T - \xi). \tag{4.6}$$

Suppose $\hat{\alpha} \in \mathcal{A}$ is a candidate for an optimal control, and denote $\hat{X}$, $(\hat{Y}, \hat{Z})$ the corresponding processes. Since the Hamiltonian is linear in $a$, we see that conditions (4.4) and (4.5) will be satisfied iff :

$$\mu_t \hat{Y}_t + \sigma_t \hat{Z}_t = 0, \quad 0 \leq t \leq T, \text{ a.s.} \tag{4.7}$$

Due to the linear-quadratic formulation of the problem, we are looking for a solution $(\hat{Y}, \hat{Z})$ to (4.6) satisfying (4.7) in the form

$$\hat{Y}_t = P_t(\hat{X}_t - Q_t), \tag{4.8}$$

where $P$ and $Q$ are adapted processes satisfying the terminal condition : $P_T = 2$ and $Q_T = \xi$. Searching for Itô processes and proceeding by identification, we see after straightforward calculation that if there exist pairs $(P, \Lambda)$ and $(Q, \Gamma)$ solutions to the BSDEs :

$$dP_t = -P_t\left[2r_t - \left(\frac{\mu_t}{\sigma_t} + \frac{\Lambda_t}{P_t}\right)^2\right] dt + \Lambda_t dW_t, \quad P_T = 2 \tag{4.9}$$

$$dQ_t = \left(r_t Q_t + \frac{\mu_t}{\sigma_t}\Gamma_t\right) + \Gamma_t dW_t, \quad Q_T = \xi, \tag{4.10}$$

then $(\hat{Y}, \hat{Z})$ is solution to the BSDE (4.6) with the mimimum condition (4.7). Existence and uniqueness of a solution to (4.9)-(4.10) are proved in [53]. Moreover, the optimal control is given by :

$$\hat{\alpha}_t = -\frac{1}{\sigma_t}\left(\frac{\mu_t}{\sigma_t} + \frac{\Lambda_t}{P_t}\right)(X_t - Q_t) + \frac{\Gamma_t}{\sigma_t}. \tag{4.11}$$

Finally, the value function is simply derived by applying Itô's formula to $P_t(X_t - Q_t)^2/2$ for $\alpha \in \mathcal{A}$ :

$$J(\alpha) = E[X_T - \xi]^2 = E\left[\frac{1}{2}P_T(X_T - Q_T)^2\right]$$

$$= \frac{1}{2}P_0(x - Q_0)^2 + \frac{1}{2}E\left[\int_0^T P_t \sigma_t^2 (\alpha_t - \hat{\alpha}_t)^2 dt\right].$$

This proves again that the optimal control is indeed given by (4.11) and the value function is equal to :

$$J(\hat{\alpha}) = \frac{1}{2}P_0(x - Q_0)^2.$$

### 4.3. Forward/backward stochastic differential equations and control

Suppose that the conditions of the sufficiency theorem 4.1 are satisfied and denote by $\hat{X}$ the corresponding optimal controlled state process, and $(\hat{Y}, \hat{Z})$ the associated adjoint dual variables. By considering the function

$$\hat{H}(t,x,y,z) \;=\; \sup_{a\in A}\mathcal{H}(t,x,a,y,z) \;=\; \sup_{a\in A}\left[b(x,a).y + \mathrm{tr}(\sigma'(x,a)z) + f(t,x,a)\right],$$

we then see, by the envelope theorem, that $(\hat{X}, \hat{Y}, \hat{Z})$ satisfy :

$$\begin{aligned}
d\hat{X}_t &= H_y(t,\hat{X}_t,\hat{Y}_t,\hat{Z}_t)dt + H_z(t,\hat{X}_t,\hat{Y}_t,\hat{Z}_t)dW_t, \\
d\hat{Y}_t &= -H_x(t,\hat{X}_t,\hat{Y}_t,\hat{Z}_t)dt + \hat{Z}_t dW_t, \quad Y_T = D_x g(\hat{X}_T).
\end{aligned}$$

where $(H_x\ H_y\ H_z)$ stands for the gradient of $H$ with respect to the $(x\ y\ z)$ variables. This is a special case of coupled forward-backward stochastic differential equation (FBSDE) written in the form :

$$\begin{aligned}
dX_t &= b(t,X_t,Y_t,Z_t)dt + \sigma(t,X_t,Y_t,Z_t)dW_t, & (4.12) \\
dY_t &= -f(t,X_t,Y_t,Z_t)dt + Z_t dW_t, \quad Y_T = h(X_T). & (4.13)
\end{aligned}$$

A solution is a triple $(X,Y,Z)$ of adapted processes with $X$ a solution to the forward s.d.e. (4.12) and $(Y,Z)$ a solution to the bsde (4.13). Coupled FBSDEs extend BSDEs since the backward component $(Y,Z)$ appear in the coefficients of the forward component $X$. The study of FBSDEs generates an important stream of research. There are essentially two methods for proving the solvability of BSDE : pure probabilistic methods like the method of continuation, see [2], [48], and combined methods of PDE and probability like the four step scheme developed in [65], or more recently extensions in [74], [27]. The connection with PDE is an extension of the Feynman-Kac formula : it states that if there exists a smooth solution with suitable Lipschitz and growth conditions to the quasilinear PDE

$$\begin{aligned}
-\frac{\partial v}{\partial t} - b(t,x,v,z(t,x,v,D_x v)).D_x v - \frac{1}{2}\mathrm{tr}(\sigma\sigma'(t,x,v,z(t,x,v,D_x v))D_x^2 v) \\
-f(t,x,v,z(t,x,v,D_x v)) = 0, \quad (t,x)\in[0,T)\times\mathbb{R}^n & (4.14) \\
v(T,x) = h(x), \quad x\in\mathbb{R}^n, & (4.15)
\end{aligned}$$

where $z(t,x,y,p)$ satisfies $z(t,x,y,p) = p'\sigma(t,x,y,z(t,x,y,p))$, then the triple $(X,Y,Z)$ determined by

$$dX_t \;=\; \tilde{b}(t,X_t)dt + \tilde{\sigma}(t,X_t)dW_t,$$

where

$$\begin{aligned}
\tilde{b}(t,x) &= b(t,x,v(t,x),z(t,x,v(t,x),D_x v(t,x))) \\
\tilde{\sigma}(t,x) &= \sigma(t,x,v(t,x),z(t,x,v(t,x),D_x v(t,x))),
\end{aligned}$$

and

$$Y_t \;=\; v(t, X_t), \qquad Z_t \;=\; z(t, X_t, v(t, X_t), D_x v(t, x_t))$$

is a solution to the FBSDE (4.12)-(4.13). Conversely, if there exists a solution $(X^{t,x}, Y^{t,x}, Z^{t,x})$ solution to the FBSDE (4.12)-(4.13) with $X_t^{t,x} = x$, then the function $v(t, x) = Y_t^{t,x}$ is a viscosity solution to (4.14)-(4.15). We refer to the lectures notes [66] for an account on the theory of FBSDE and their applications.

## 5. Numerical issues

Explicit analytic solutions to stochastic control problems are rare and one has to resort to numerical approximations. This has led to an important area of research on numerical methods for stochastic optimization. Motivated in particular by problems in quantitative finance, one faces challenging numerical problems arising in portfolio optimization in high dimension, under partial information, under transaction costs ... We distinguish two types of numerical methods for stochastic control problems : purely deterministic and probabilistic methods. We briefly survey some of the developments and advances in these fields.

### *5.1. Deterministic methods*

Purely deterministic methods from numerical analysis consist in discretizing the nonlinear HJB PDE satisfied by the value function of the stochastic control problem. The discretization by finite difference methods or finite elements methods leads to an approximation of the value function at the points of the space-time grid. Convergence of the numerical scheme is proved rigorously by stability results for viscosity solutions, see [7]. There are also recent results on convergence rate for the finite difference approximation of the HJB equation in [57] and [6]. Computational methods are studied e.g. in [1]. For some illustrations in financial problems, see e.g. [35] or [88]. From a computational viewpoint, the limitation of purely deterministic methods is the dimension of the state variable and in practice, calculations are done for low dimensions, say 1 or 2.

### *5.2. Probabilistic methods*

Probabilistic numerical methods for optimization problems considerably developed over these years. A classical method, based on the dynamic programming principle, is the Markov chain approximation introduced by Kushner [58], see also the book [59] for an account of this method. This consists in approximating the original continuous-time controlled process by a suitable controlled Markov chain on a lattice satisfying a local consistency condition. The stochastic control problem is then numerically solved from the dynamic programming principle applied to the Markov chain approximation, which leads to a backward recursive algorihm. The finite difference scheme is a typical example of a numerical scheme for an approximating Markov chain with nearest neighbor transitions. In this method, the required stability condition may be very restrictive in the case of control-dependent diffusion. On the other hand, the lattice is regardless of the structure of the Markov chain. Moreover, its size is growing exponentially with the dimension. So, although the Markov chain is easily implemented, calculations can be done in practice only for low dimensions as in the purely deterministic approach. To overcome the dimensionality problem, Pagès, Pham and Printems [71] propose a probabilistic method based on optimal quantization methods for numerically solving stochastic control problems in dimensions larger than 3. Like in usual probabilistic methods, we start from a time discretization of the controlled problem : the process $(X_t)$ is approximated by its Euler scheme denoted $\bar{X}_k$ at time $t_k = kT/n$. Then, and in the spirit of the Markov chain approximation method, the Euler scheme is approximated at every date $k$ by a process $\hat{X}_k$, taking finitely many states. The optimal quantization approach consists in finding the best approximation according to $L^p$-norm of $\bar{X}_k$ by $\hat{X}_k$. This is achieved by stochastic gradient descent method based on Monte-Carlo simulations of $\bar{X}_k$, which also allows to estimate the probability transitions of the grid points. The control problem is then numerically implemented from the backward dynamic programming formula. The main interest in optimal quantization is that given a total number of points to be dispatched among all grids, one gets optimal grids with respect to the distribution structure of the original process. This method is numerically tested in dimension 3 for the mean-variance hedging problem with stochastic volatility. We refer to [72] for an account of other applications of quantization methods to numerical problems in finance.

Another approach is to solve numerically the FBSDE associated via the maximum principle to the stochastic control problem. Several approximation methods to FBSDEs were proposed in the literature. We mention for BSDE the Markov chain approximation in [64], regression methods in [43], quantization methods in [4], the Monte-Carlo Malliavin method in [18], and for FBSDE the four-step scheme in [29], and quantization method in [28].

## 6. Conclusion

Many theoretical and numerical advances have been recently realized in the field of stochastic control. They also contribute to the fields of nonlinear partial differential equations and backward stochastic differential equations. Present developments include modelling with jump-diffusion processes, which lead to integrodifferential equations, and FBSDEs with reflecting barriers arising typically in optimal stopping problems. From a numerical viewpoint, challenging problems are the search for fast and efficient schemes for stochastic control problem in high dimension arising in quantitative finance. In this direction, partial observation problems, which lead to infinite-dimensional problem, are very few studied numerically and represent an important stream of research. Classical control problems in finance are optimal portfolio allocation but there is a large potential applicability for other economics and finance problems such as in contract or game theory, credit risk or liquidity risk models. In return, these questions should raise new developments in terms of mathematical theories, still in accordance with meaningful applied problems.

## References

[1] Akian M. (1990) : "Analyse de l'algorithme multigrille FMGH de résolution d'équations d'Hamilton-Jacobi-Bellman", *Analysis and Optimization of systems*, Lect. Notes in Contr. and Inf. Sciences, **144**, Springer-Verlag, pp. 113-122.

[2] Antonelli F. (1993) : "Backward-forward stochastic differential equations", *Annals of Appl. Prob.*, **3**, 777-793.

[3] Artzner P., Delbaen F., Eber J.M. and D. Heath (1999) : "Coherent measures of risk", *Mathematical Finance*, **9**, 203-228.

[4] Bally V. and G. Pagès (2003) : "Error analysis of the optimal quantization algorithm for obstacle problems", *Stochastic Processes and their Applications*, **106**, 1-40.

[5] Baras J., Elliott R. and M. Kohlmann (1989) : "The partially obsered stochastic minimum principle", *SIAM J. Cont. Optim.*, **27**, 1279-1292.

[6] Barles G. and E. Jakobsen (2004) : "Error bounds for monotone approximations schemes for Hamilton-Jacobi-Bellman equations", to appear in *SIAM J. Num. Anal.*.

[7] Barles G. and P. Souganidis (1991) : : "Convergence of approximation schemes for fully non linear second-order equations", *Asymptotics Analysis*, **4**, pp.271-283.

[8] Bellman R. (1957) : Dynamic programming, Princeton university press.

[9] Bensoussan A. and J.L. Lions (1982) : Contrôle impulsionnel et inéquations variationnelles, Dunod.

[10] Bensoussan A. (1992) : Stochastic control of partially observable systems, Cambridge University Press.

[11] Bensoussan A. and H. Nagai (1991) : "An ergodic control problem arising from the principal eigenfunction of an elliptic operator", *J. Math. Soc. Japan*, **43**, 49-65.

[12] Bielecki T. and S. Pliska (1999) : "Risk-sensitive dynamic asset management", *Applied Math. Optim.*, **39**, 337-360. MR1675114

[13] Blanchet-Scalliet C., El Karoui N., Jeanblanc M. and L. Martellini (2002) : "Optimal investment and consumption when time-horizon is uncertain", Preprint.

[14] Bobrovnytska O. and M. Schweizer (2004) : "Mean-variance hedging and stochastic control : beyond the Brownian setting", *IEE Trans. on Aut. Cont.*, **49**, 1-14.

[15] Borkar V. (1989) : Optimal control of diffusion processes, Pitman Research Notes in Math., 203. MR1005532

[16] Borkar V. (2005) : "Controlled diffusion processes", *Probability surveys*, **2**, 213-244.

[17] Bouchard B. and H. Pham (2004) : "Wealth-path dependent utility maximization in incomplete markets", *Finance and Stochastics*, **8**, 579-603.

[18] Bouchard B. and N. Touzi (2004) : "Discrete time approximation and Monte Carlo simulation of backward stochastic differential equations", *Stochastic Processes and their Applications*, **111**, 175-206.

[19] Brekke K. and B. Oksendal (1994) : "Optimal switching in an economic activity under uncertainty", *SIAM J. Cont. Optim.*, **32**, 1021-1036.

[20] Broadie M., Cvitanic J. and M. Soner (1998) : "Optimal replication of contingent claims under portfolio constraints", *Review of Fin. Studies*, **11**, 59-79.

[21] Carmona R. and N. Touzi (2004) : "Optimal multiple stopping and the valuation of swing options", to appear in *Mathematical Finance*.

[22] Çetin U., Jarrow R. and P. Protter (2004) : "Liquidity risk and arbitrage pricing theory", *Finance and Stochastics*, **8**, 311-341.

[23] Crandall M., Ishii. H and P.L. Lions (1992) : "User's Guide to Viscosity Solutions of Second Order Partial Differential Equations", *Bull. Amer. Math. Soc.*, **27**, 1-67.

[24] Cvitanic J., Pham H. and N. Touzi (1999a) : "Superreplication in stochastic volatility models under portfolio constraints", *Journal of Appplied Probability*, **36**, 523-545. MR1724796

[25] Cvitanic J., Pham H. and N. Touzi (1999b) : "A closed form solution for the super-replication problem under transaction costs", *Finance and Stochastics*, **3**, 35-54.

[26] Davis M. and P. Varaiya (1973) : "Dynamic programming conditions for partially observable systems", *SIAM J. Cont.*, **11**, 226-261.

[27] Delarue, F. (2002) : "On the existence and uniqueness of solutions to FBSDEs in a non-degenerate case". *Stoch. Process. Appl.*, **99**, 209-286.

[28] Delarue F. and S. Menozzi (2004) : "A forward-backward stochastic algorithm for quasi-linear PDEs", to appear in *Annals of Applied Probability*.

[29] Douglas J., Ma J. and P. Protter (1996) : "Numerical methods for forward-backward stochastic differential equations", *Annals of Applied Probability*, **6**, 940-968.

[30] Duckworth K. and M. Zervos (2001) : "A model for investment decisions with switching costs", *Annals of Applied Probability*, **11**, 239-250.

[31] El Karoui N. (1981) : Les Aspects Probabilistes du Contrôle Stochastique, Lect. Notes in Math., 816, Springer Verlag.

[32] El Karoui N. and L. Mazliak (editors) (1997) : Backward stochastic differential equations, Pitman research notes in mathematics series.

[33] El Karoui N., S. Peng and M.C. Quenez (1997) : "Backward stochastic differential equations in finance", *Mathematical Finance*, **7**, 1-71. MR1434407

[34] El Karoui N. and M.C. Quenez (1995) : "Dynamic programming and pricing contingent claims in incomplete markets", *SIAM J. Cont. Optim.*, **33**, 29-

66. MR1311659
[35] Fitzpatrick B. and W. Fleming (1991) : "Numerical methods for an optimal Investment-Consumption Model", *Mathematics of Operation Research*, **16**, pp.823-841.
[36] Fleming W. (1968) : "Optimal control of partially observable systems", *SIAM J. Cont. Optim.*, **6**, 194-214.
[37] Fleming, W. and D. Hernandez-Hernandez (2003) : "An optimal consumption model with stochastic volatility", *Finance Stoch.*, **7**, 245-262.
[38] Fleming W. and W. McEneaney (1995) : "Risk-sensitive control on an infinite horizon", *SIAM J. Cont. and Optim.*, **33**, 1881-1915.
[39] Fleming W. and R. Rishel (1975) : Deterministic and stochastic optimal control, Springer Verlag.
[40] Fleming W. and S. Sheu (2000) : "Risk sensitive control and an optimal investment model", *Math. Finance*, **10**, 197-213. MR1802598
[41] Fleming W. and M. Soner (1993) : Controlled Markov processes and viscosity solutions, Springer Verlag.
[42] Fuhrman M., Hu Y. and G. Tessitore (2005) : "On a class of stochastic optimal control problems related to BSDEs with quadratic growth", Preprint.
[43] Gobet E., Lemor J.P. and X. Warin (2004) : "A regression-based Monte-Carlo method for backward stochastic differential equations", to appear in *Annals of Applied Probability*.
[44] Gundel A. (2004) : "Robust utility maximization for complete and incomplete market models", *Finance and Stochastics*, **9**, 151-176.
[45] Guo X. (2001) : "An explicit solution to an optimal stopping problem with regime switching", *Journal of Applied Probability*, **38**, 464-481.
[46] Guo X. and H. Pham (2005) : "Optimal partially reversible investment with entry decision and general production function", *Stoc. Proc. Appli.*, **115**, 705-736. MR2132595
[47] Hata H. and J. Sekine (2005) : "Solving a large deviations control problem with a nonlinear factor model", Preprint.
[48] Hu Y. and S. Peng (1995) : "Solution of forward-backward stochastic differential equations", *Prob. Theo. Rel. Fields*, **103**, 273-283.
[49] Jeanblanc M. and A. Shiryaev (1995) : "Optimization of the flow of dividends", *Russian Math Surveys*, **50**, 257-277.
[50] Kabanov, Yu. and C. Kluppelberg (2004) : "A geometric approach to portfolio optimization in models with transaction costs", *Finance and Stochastics*, **8**, 207-227.
[51] Karatzas I. (1980) : "On a stochastic representation for the principal eigenvalue of a second order differential equation", *Stochastics and Stochastics Reports*, **3**, 305-321.
[52] Kobylanski M. (2000) : "Backward stochastic differential equations and partial differential equations with quadratic growth", *Annals of Probability*, **28**, 558-602.
[53] Kohlmann M. and S. Tang (2002) : "Global adapted solution of one-dimensional backward stochastic differential Riccati equations, with application to the mean-variance hedging", *Stochastic Process. Appl.*, **97**, 255–

288.
[54] Kohlmann M. and X.Y. Zhou (2000) : "Relationship between backward stochastic differential equations and stochastic controls : a linear-quadratic approach", *SIAM Journal on Control and Optimization*, **38**, 1392-1407.
[55] Korn R. (1998) : "Portfolio optimization with strictly positive transaction costs and impulse control", *Finance and Stochastics*, **2**, 85-114.
[56] Krylov N. (1980) : Controlled Diffusion Processes, Springer Verlag.
[57] Krylov N. (2000) : "On the rate of convergence of finite difference approximations for Bellman's equations with variable coefficients", *Probab. Theory Relat. Fields*, **117**, 1-16.
[58] Kushner H.J. (1977) : "Approximation and weak convergence methods for random processes, with applications to stochastic systems theory", MIT Press Series in Signal Processing, Optimization, and Control, **6**, MIT Press, Cambridge, MA, 1984, 269 pp.
[59] Kushner H.J. and P. Dupuis (2001) : Numerical methods for stochastic control problems in continuous time, $2^{nd}$ edition, Applications of Mathematics, **24**, Stochastic Modelling and Applied Probability, Springer-Verlag, New York, 475 pp.
[60] Ladyzhenskaya O., Solonnikov V. and N. Uralseva (1968) : Linear and quasilinear equations of parabolic type, American Mathematical Society, Providence.
[61] Lions P.L. (1983) : "Optimal control of diffusion processes and Hamilton-Jacobi-Bellman equations"", *Comm. P.D.E*, **8**, Part I, 1101-1134, Part II, 1229-1276
[62] Lions P.L. (1988) : "Viscosity solutions and optimal stochastic control in infinite dimension", Part I *Acta Math.*, **161**, 243-278.
[63] Ly Vath V., Mnif M. and H. Pham (2005) : "A model of optimal portfolio selection under liquidity risk and price impact", preprint PMA, University Paris 6-Paris 7.
[64] Ma J., Protter P., San Martin J. and S. Torres (2002) : "Numerical method for backward stochastic differential equations", *Ann. Appl. Prob.*, **12**, 302-316.
[65] Ma J., Protter P. and J. Yong (1994) : "Solving forward-backward stochastic differential equations explicitly : a four step scheme", *Prob. Theo. Rel. Fields*, **98**, 339-359.
[66] Ma J. and J. Yong (2000) : Forward-backward stochastic differential equations and their applications, Lect. Notes in Math., 1702.
[67] Mania M. and R. Tevzadze (2003) : "Backward stochastic PDE and imperfect hedging", *Int. J. Theo. App. Fin.*, **6**, 663-692.
[68] Nisio M. (1981) : Lectures on Stochastic Control Theory, ISI Lect. Notes, 9, Kaigai Publ. Osaka.
[69] Oksendal B. and A. Sulem (2002) : "Optimal consumption and portfolio with both fixed and proportional transaction costs", *SIAM J. Cont. Optim.*, **40**, 1765-1790.
[70] Oksendal B. and A. Sulem (2004) : Applied stochastic control of jump diffusion, Springer Verlag. MR2109687

[71] Pagès G., Pham H. and J. Printems (2004a) : "An optimal Markovian quantization algorithm for multi-dimensional stochastic control problems", *Stochastics and Dynamics*, **4**, 501-545.

[72] Pagès G., Pham H. and J. Printems (2004b) : "Optimal quantization methods and applications to numerical problems in finance", in Handbook of computational and numerical methods in finance, ed. S. Rachev, Birkhäuser. MR2083048

[73] Pardoux E. and S. Peng (1990) : "Adapted solutions of a backward stochastic differential equation", *Systems and Control Letters*, **14**, 55-61.

[74] Pardoux E. and S. Tang (1999) : "Forward-backward stochastic differential equations and quasilinear parabolic pdes ", *Prob. Theo. Rel. Fields.*, **114**, 123-150.

[75] Peng S. (1990) : "A general stochastic maximum principle for optimal control diffusions", *SIAM J. Cont. Optim.*, **28**, 966-979.

[76] Pham H. (2000) : "On quadratic hedging in continuous time", *Math. Meth. Oper. Res.*, **51**, 315-339.

[77] Pham H. (2002) : "Smooth solutions to optimal investment models with stochastic volatilities and portfolio constraints", *Appl. Math. Optim.*, **46**, 55-78.

[78] Pham H. (2003a) : "A large deviations approach to optimal long term investment'", *Finance and Stochastics*, **7**, 169-195.

[79] Pham H. (2003b) : "A risk-sensitive control dual approach to a large deviations control problem", *Systems and Control Letters*, **49**, 295-309.

[80] Pham H. (2005a) : "On the smooth-fit property for one-dimensional optimal switching problem", to appear in *Séminaire de Probabilités*, Vol. XL.

[81] Pham H. (2005b) : Optimisation et contrôle stochastique appliqués à la finance, to appear, Series Mathématiques et Applications, Springer Verlag.

[82] Schied A. (2003) : "Optimal investment for robust utility functionals in complete markets", to appear in *Mathematics of Operations Research*.

[83] Schweizer M. (2001) : "A guided tour through quadratic hedging approaches", in *Option pricing, Interest rates and risk management*, J. Cvitanic, E. Jouini, and m. Musiela eds., Cambridge Univ. Press.

[84] Shreve S. and M. Soner (1994) : "Optimal investment and consumption with transaction costs", *Annals of Applied Probability*, **4**, 609-692.

[85] Soner M. and N. Touzi (2000) : "Super replication under gamma constraints", SIAM Journal on Control and Optimization, **39**, 73-96.

[86] Soner M. and N. Touzi (2002) : "Stochastic target problems, dynamic programming and viscosity solutions", SIAM Journal on Control and Optimization, **41**, 404-424.

[87] Tang S. and J. Yong (1993) : "Finite horizon stochastic optimal switching and impulse controls with a viscosity solution approach", *Stoch. and Stoch. Reports*, **45**, 145-176. MR1306930

[88] Tourin A. and T. Zariphopoulou (1994) : "Numerical schemes for investment models with singular transactions", *Computational Economics*, **7**, pp.287-307.

[89] Yong J. and X.Y. Zhou (2000) : Stochastic controls, Hamiltonian systems

and HJB equations, Springer Verlag. MR1696772
[90] Zariphopoulou T. (1988) : Optimal investment-consumption models with constraints, Brown University, Phd.
[91] Zariphopoulou T. (2001) : "A solution approach to valuation with unhedgeable risks", *Finance and Stochastics*, **5**, 61-82. MR1807876
[92] Zitkovic G. (2005) : "Utility Maximization with a Stochastic Clock and an Unbounded Random Endowment", *Annals of Applied Probability*, **15**, 748-777.
[93] Zhou X.Y. (1993) : "On the necessary conditions of optimal controls for stochastic partial differential equations", *SIAM J. Cont. Optim.*, **31**, 1462-1478.